\documentclass[final,onefignum,onetabnum,onealgnum]{siamart251216}

\usepackage{amsmath,amsfonts,amsopn,amssymb}
\usepackage{graphicx}
\usepackage[caption=false]{subfig}
\usepackage{multirow,relsize,makecell}
\usepackage{url}
\ifpdf
  \DeclareGraphicsExtensions{.eps,.pdf,.png,.jpg}
\else
  \DeclareGraphicsExtensions{.eps}
\fi

\usepackage{xr}

\usepackage{algorithm,algorithmic}
\usepackage{accents}
\usepackage{tikz}
\usetikzlibrary{arrows.meta,positioning,backgrounds,overlay-beamer-styles,calc}

\numberwithin{theorem}{section}
\newsiamremark{remark}{Remark}
\newsiamremark{assumption}{Assumption}
\newsiamremark{example}{Example}
\crefname{remark}{Remark}{Remarks}
\crefname{assumption}{Assumption}{Assumptions}
\crefname{example}{Example}{Examples}

\headers{Subspace correction for convex optimization}{Boou Jiang, Jongho Park, and Jinchao Xu}

\title{Subspace correction as a general framework\\ for convex optimization algorithms\thanks{Submitted to arXiv.\funding{This work was supported in part by the KAUST Baseline Research Fund (BAS/1/1699)
and the University Development Fund from The Chinese University of Hong Kong,
Shenzhen (UDF01004630).}}
}

\author{
Boou Jiang\thanks{Applied Mathematics and Computational Sciences Program, Computer, Electrical and Mathematical Science and Engineering Division, King Abdullah University of Science and Technology~(KAUST), Thuwal 23955, Saudi Arabia
(\email{boou.jiang@kaust.edu.sa}, \email{jinchao.xu@kaust.edu.sa}).}
\and
Jongho Park\footnotemark[2]\hspace{0.35em}\thanks{School of Science and Engineering, The Chinese University of Hong Kong, Shenzhen, Guangdong 518172, China
 (\email{jonghopark@cuhk.edu.cn}).}
\and
Jinchao Xu\footnotemark[2]
 }

\ifpdf
\hypersetup{
  pdftitle={Subspace correction as a general framework for convex optimization algorithms},
  pdfauthor={Boou Jiang, Jongho Park, and Jinchao Xu}
}
\fi

\usepackage{color}

\begin{document}
\maketitle

\begin{abstract}
This paper shows that subspace correction methods provide a common algorithmic and theoretical foundation for several classes of convex optimization algorithms, including operator splitting, alternating projection, and multiplier methods.
The underlying principle is to decompose a problem into smaller subproblems and combine their solutions, a strategy that appears throughout iterative algorithms.

The main tool is an iterate-level formalism, which we call dualization, that abstracts classical primal--dual correspondences and relates a subspace correction method for a dual problem to an algorithm for the corresponding primal problem via a primal--dual consistency relation. 
At the algorithmic level, dualizing successive subspace correction yields the Peaceman--Rachford and Douglas--Rachford splitting methods; the von Neumann and Dykstra alternating projection algorithms arise as special cases. 
Dualizing parallel subspace correction yields a parallel splitting method. Multiplier methods, including the alternating direction method of multipliers~(ADMM), are connected to subspace correction through the same mechanism together with equivalent block formulations. 
In particular, a multi-block ADMM-type algorithm is obtained by dualizing a splitting method derived from successive subspace correction.

At the theoretical level, the primal--dual consistency relation transfers convergence estimates from subspace correction to the derived algorithms under suitable assumptions. 
Thus, these operator splitting, alternating projection, and multiplier algorithms can be derived from subspace correction at both the algorithmic and theoretical levels.
\end{abstract}

\begin{keywords}
Subspace correction method,
Operator splitting method,
Alternating projection method,
ADMM,
Dualization,
Convex duality,
Convex optimization 
\end{keywords}

\begin{AMS}
90C46,  
49M27,  
68W15,  
65N55, 
90C25   
\end{AMS}

\section{Introduction}
\label{Sec:Introduction}
Convex optimization plays a fundamental role in various fields of science and engineering, including scientific computing, machine learning, and operations research.
For recent advances in numerical algorithms for convex optimization and related problems, see, e.g.,~\cite{CP:2016a,CKCH:2023,DST:2021,Nesterov:2018}.

Because convex optimization arises in such diverse fields, several algorithms have been developed independently and only later recognized as equivalent.
The proximal point algorithm~\cite{Martinet:1970,Rockafellar:1976}, the augmented Lagrangian method~\cite{Hestenes:1969,Powell:1969}, and the split Bregman method~\cite{COS:2010,GO:2009} were introduced for solving monotone inclusions, constrained optimization problems, and regularized problems, respectively. Despite these different motivations, they share a common structure~\cite{EB:1992,Setzer:2011,WT:2010}.
Primal--dual algorithms~\cite{CP:2011,EZC:2010} solve saddle point problems, the alternating direction method of multipliers~(ADMM)~\cite{BPCPE:2011,Glowinski:2014} treats linearly constrained problems by alternating minimization, and Douglas--Rachford splitting~\cite{DR:1956,LM:1979} separates problems involving sums of operators. Relations among these methods have been studied in~\cite{CP:2011,CKCH:2023,GL:1989}.
These examples motivate the search for a common framework that explains both the construction of the algorithms and their convergence theories.

A similar question arises in scientific computing, particularly for iterative methods for linear systems obtained by discretizing partial differential equations~(PDEs).
A central framework is provided by subspace correction methods~\cite{Xu:1992}, also called Schwarz methods in the domain decomposition literature~\cite{TW:2005}.
The solution space is decomposed into a sum of subspaces, local problems are solved on these subspaces, and the local solutions are combined to update the global iterate.
This framework includes classical Jacobi and Gauss--Seidel iterations, multigrid methods~\cite{Bramble:1993,BPX:1990}, and domain decomposition methods~\cite{TW:2005,XZ:1998}.
A sharp convergence characterization is given by the so-called Xu--Zikatanov identity~\cite{Brenner:2013,LWXZ:2008,XZ:2002}.
Subspace correction has been extended from linear systems to smooth and nonsmooth convex optimization~\cite{Carstensen:1997,TE:1998,TX:2002}; more recent developments include~\cite{Park:2020,Park:2022a}.
Applications include nonlinear PDEs~\cite{CHW:2020,Park:2024a} and total-variation minimization~\cite{CTWY:2015,HL:2025,LP:2020}.

The central viewpoint of this paper is that subspace correction provides a common framework not only for algorithms in scientific computing but also for several classes of convex optimization algorithms, including operator splitting, alternating projection, and multiplier methods.

We will use the following three model problems to discuss various optimization algorithms, including subspace correction, operator splitting, and multiplier methods:
\begin{enumerate}
\item
The general convex minimization problem
\begin{equation}
\label{model_problem_1}
\min_{v\in V}E(v).
\end{equation}
Here, $V$ is a real Hilbert space and $E$ is a convex function on $V$.

\item
The structured composite convex optimization problem
\begin{equation}
\label{model_problem_2}
\min_{v\in V}
\left\{
F(v)+G(Bv)
\right\}.
\end{equation}
Here, $V$ and $W$ are real Hilbert spaces, $F$ and $G$ are convex functions on $V$ and $W$, respectively, and $B\colon V\to W$ is a bounded linear operator.
Clearly,~\eqref{model_problem_2} is a special case of~\eqref{model_problem_1}.

\item
The linearly constrained convex optimization problem
\begin{equation}
\label{model_problem_3}
\min_{v\in V}F(v)
\quad\text{subject to}\quad
Bv=g.
\end{equation}
Here, $V$ and $W$ are real Hilbert spaces, $F$ is a convex function on $V$, $B\colon V\to W$ is a bounded linear operator, and $g\in W$.
Unlike the first two problems,~\eqref{model_problem_3} is not written directly in the unconstrained form~\eqref{model_problem_1}.
We will see later that, through convex duality, it can nevertheless be treated within the same framework.
\end{enumerate}
The precise assumptions used for these model problems are stated when they are needed.

Subspace correction methods were originally developed for linear problems~\cite{Xu:1992,XZ:2002} and have since been extended to convex optimization problems of the form~\eqref{model_problem_1}~\cite{Carstensen:1997,TE:1998,TX:2002}.
They include parallel subspace correction~(PSC), in which all local problems are solved from the same iterate, and successive subspace correction~(SSC), in which local corrections are applied sequentially.
On a suitable product space, PSC and SSC are equivalent to block Jacobi and block Gauss--Seidel methods, respectively, for an \emph{expanded formulation}~\cite{PX:2025}.
While block-coordinate methods are typically formulated in terms of disjoint variable blocks, subspace decompositions may overlap or be organized hierarchically, as in domain decomposition and multigrid methods.
The expanded formulation provides a connection between such subspace structures and block optimization methods.

In~\eqref{model_problem_2}, when the term $G\circ B$ is given by a sum of multiple terms, operator splitting methods are commonly used to decompose the problem so that the individual terms can be handled separately.
Operator splitting originated from the Peaceman--Rachford and Douglas--Rachford methods for PDEs~\cite{DR:1956,PR:1955}; see also~\cite{BVY:1962}.
It has since become a standard tool for convex optimization and monotone operator problems~\cite{EB:1992,LM:1979}.
The broader framework includes forward--backward algorithms~\cite{BT:2009,CP:2016a} and primal--dual algorithms~\cite{CP:2011,CP:2016b}; see also~\cite{CP:2009,CKCH:2023}.

Alternating projection methods arise as special cases when the terms in $G\circ B$ are indicator functions of convex sets.
They repeatedly project onto the individual sets to solve feasibility or best-approximation problems over their intersection.
The von Neumann algorithm alternates projections between closed subspaces~\cite{Neumann:1950}, whereas the Dykstra algorithm introduces correction terms to compute the projection onto the intersection of general closed convex sets~\cite{Dykstra:1983}.
Connections between these methods and subspace correction or block coordinate descent have been studied in~\cite{GM:1989,Han:1988,Tibshirani:2017}; see also~\cite{XZ:2002}.
A comprehensive treatment of alternating projection methods can be found in~\cite{BB:1996}.

Problem~\eqref{model_problem_3} can be treated by multiplier methods, which introduce Lagrange multipliers to enforce the linear constraint.
When the primal variable has a block structure, the resulting subproblems can be treated blockwise, leading naturally to ADMM-type methods.
The augmented Lagrangian method adds a quadratic penalty term to the Lagrangian~\cite{Hestenes:1969,Powell:1969}, while decomposing the primal update into blockwise subproblems leads to ADMM~\cite{GM:1976,GM:1975}.
For general reviews of ADMM, see~\cite{BPCPE:2011,Glowinski:2014}.

At first sight, subspace correction, operator splitting, and multiplier methods appear to constitute separate algorithmic families, with alternating projection arising as a special case of operator splitting. 
The main purpose of this paper is to show that they instead form a hierarchy with subspace correction at its foundation.
To establish this hierarchy, we formalize an iterate-level mechanism, called \emph{dualization}, which associates a subspace correction method for a dual problem with an algorithm for the corresponding primal problem.
The corresponding iterates satisfy a primal--dual consistency relation obtained from Fenchel--Rockafellar duality~\cite{Rockafellar:1967,Rockafellar:1970}.

The term dualization is not intended to suggest that the underlying primal--dual transformations are new.
Rather, it abstracts a mechanism already implicit in classical correspondences, including the links between Douglas--Rachford splitting and ADMM~\cite{EB:1992,Setzer:2011}, and among the proximal point algorithm, the augmented Lagrangian method, and the split Bregman method~\cite{COS:2010,GO:2009,WT:2010}.
Related duality-based constructions have also appeared in specific applications~\cite{LG:2019,LN:2017,PX:2023}, including the dualized level-set method in~\cite{MMP:2025}.

Algorithmically, dualizing SSC yields Peaceman--Rachford and Douglas--Rachford splitting methods, while dualizing PSC yields a parallel splitting method.
Alternating projection methods arise as special cases, and the multiplier methods are obtained through the same construction or through a chain of dualizations and equivalent block formulations.
Theoretically, the iterate-level consistency relation transfers convergence estimates for subspace correction methods to the derived algorithms under suitable assumptions.
Thus, our main conclusion is that subspace correction generates the operator splitting, alternating projection, and multiplier algorithms described above and, under suitable assumptions, provides their convergence theory.

The principal results supporting this conclusion are summarized as follows:
\begin{enumerate}
\item
\Cref{Thm:convergence_transfer} provides the theoretical mechanism for transferring convergence results.
Under suitable conditions, convergence of a dual sequence implies convergence of its dualization, together with a corresponding convergence estimate.
Consequently, convergence results for subspace correction methods can be transferred to their dualizations under the stated assumptions.

\item
\Cref{Thm:expanded_system} identifies subspace correction methods with block methods on an expanded product space.
This equivalence provides the basic connection between subspace decompositions and block algorithms.

\item
\Cref{Thm:multiple_convex_PR,Thm:multiple_convex_DR} show that the dualization-based Peaceman--Rachford and Douglas--Rachford methods considered in \cref{Sec:Splitting} are obtained by dualizing the corresponding SSC methods.
Their alternating-projection specializations recover the relations between the von Neumann algorithm and SSC~\cite{XZ:2002}, and between the Dykstra algorithm and block coordinate descent~\cite{GM:1989,Han:1988}.
Similarly, \cref{Thm:multiple_convex_parallel} shows that the parallel splitting method is obtained by dualizing PSC.

\item
\Cref{Thm:ADMM_multiple} shows that the ADMM considered in \cref{Sec:ADMM} is obtained by dualizing the corresponding Peaceman--Rachford method.
\end{enumerate}

The rest of this paper is organized as follows.
In \cref{Sec:Convex}, we provide preliminary notation and a concise overview of convex analysis.
In \cref{Sec:Dualization}, we review Fenchel--Rockafellar duality, introduce dualization, and present several motivating applications.
In \cref{Sec:MSC}, we review subspace correction methods for linear systems, their interpretation as block methods for an expanded system, and their extension to convex optimization.
In \cref{Sec:Splitting}, we derive sequential, relaxed, and parallel splitting methods from subspace correction and treat alternating projection methods as special cases.
In \cref{Sec:ADMM}, we derive an ADMM-type method from an operator splitting method via dualization.
Finally, in \cref{Sec:Conclusion}, we conclude with some remarks.

\section{Preliminaries}
\label{Sec:Convex}
In this section, we introduce notation for linear operators and briefly review essential concepts from convex analysis.
For a more comprehensive treatment of convex analysis relevant to this work, we refer the reader to the monographs~\cite{BC:2011,Rockafellar:1970,RW:2009}.

\subsection{Linear operators}

Throughout this paper, $V$, $W$, and all other abstract spaces are real Hilbert spaces, and $V^*$ and $W^*$ denote their dual spaces.
By the Riesz representation theorem, every element of the dual of a real Hilbert space has a unique Riesz representative in that Hilbert space.
We use the same symbol for a dual element and its Riesz representative, while retaining the notation $V$ and $V^*$, and similarly for the other spaces, to distinguish their primal and dual roles.
All linear operators between these spaces are assumed to be bounded.
We use the same notation $(\cdot,\cdot)$ for inner products and duality pairings, and the same notation $\|\cdot\|$ for norms, with the underlying spaces distinguishing their respective meanings~(cf.~\cite{XZ:2017}).
Thus, for $q\in V^*$ and $v\in V$, the action of $q$ on $v$ is denoted by $(q,v)$.
An operator $A\colon V\to V^*$ is called symmetric positive definite~(SPD) if it is symmetric and there exists $c_A>0$ such that
\begin{equation*}
(Av,v)\geq c_A\|v\|^2,
\quad v\in V.
\end{equation*}
The analogous convention is used for operators from $V^*$ to $V$.
For a linear operator $A\colon V\to W$, we denote its adjoint by $A^*\colon W^*\to V^*$, i.e.,
\begin{equation*}
    (q, Av) = (A^* q, v),
    \quad v \in V,\ q \in W^*.
\end{equation*}

\subsection{Convex functions}

We write $\overline{\mathbb{R}}=\mathbb{R}\cup\{\pm\infty\}$ for the extended real line.
Let $F\colon V\to\overline{\mathbb{R}}$ be an extended real-valued function.
The \emph{effective domain} of $F$, denoted by $\operatorname{dom}F$, is defined as
\begin{equation*}
\operatorname{dom}F
=
\{v\in V:F(v)<\infty\}.
\end{equation*}
For example, given a subset $K\subset V$, its \emph{indicator function}
$\chi_K\colon V\to\overline{\mathbb{R}}$ is defined by
\begin{equation*}
\chi_K(v)
=
\begin{cases}
0, & \text{if }v\in K,\\
\infty, & \text{if }v\notin K,
\end{cases}
\quad v\in V,
\end{equation*}
and satisfies $\operatorname{dom}\chi_K=K$.
The function $F$ is said to be \emph{proper} if
$\operatorname{dom}F\neq\emptyset$ and $F$ never takes the value
$-\infty$.
It is said to be \emph{lower semicontinuous} if, for every $v\in V$ and every sequence $\{v_n\}$ converging to $v$, we have
$F(v)\leq\liminf_{n\to\infty}F(v_n)$.

An extended real-valued function
$F\colon V\to\overline{\mathbb{R}}$ is said to be \emph{convex} if
\begin{equation*}
F((1-t)v+tw)
\leq
(1-t)F(v)+tF(w),
\quad
v,w\in V,\ t\in[0,1],
\end{equation*}
with the usual conventions of extended-real arithmetic~(see, e.g.,~\cite{RW:2009}).
If $F\colon V\to\mathbb{R}$ is Fr\'echet differentiable, then it is convex if and only if
\begin{equation}
\label{first-order-convexity}
F(v)
\geq
F(w)+(F'(w),v-w),
\quad
v,w\in V,
\end{equation}
where $F'(w)\in V^*$ denotes the derivative of $F$ at $w$.
Thus, the graph of a differentiable convex function lies above each of its tangent affine functionals.

A Fr\'echet differentiable convex function
$F\colon V\to \mathbb{R}$ is said to be
$\mu$-\emph{strongly convex}, for some $\mu>0$, if 
\begin{equation}
\label{strong_convexity}
F(v)
\geq
F(w)+(F'(w),v-w)
+
\frac{\mu}{2}\|v-w\|^2,
\quad v, w \in V.
\end{equation}
On the other hand, it is said to be
$L$-\emph{smooth}, for some $L>0$, if its derivative is
$L$-Lipschitz continuous as a map from $V$ to $V^*$, i.e.,
\begin{equation*}
\|F'(v)-F'(w)\|
\leq
L\|v-w\|,
\quad
v,w\in V.
\end{equation*}

Finally, given a Fr\'echet differentiable convex function
$F\colon V\to\mathbb{R}$, the \emph{Bregman divergence} $D_F$
associated with $F$ is defined as
\begin{equation*}
D_F(v;w)
=
F(v)-F(w)-(F'(w),v-w),
\quad
v,w\in V.
\end{equation*}

The first-order characterization~\eqref{first-order-convexity} extends to nondifferentiable convex functions through the notion of a subdifferential.
For a convex function $F\colon V\to\overline{\mathbb{R}}$ and $w\in\operatorname{dom}F$, the
\emph{subdifferential} of $F$ at $w$ is defined as
\begin{equation*}
\begin{aligned}
\partial F(w)
&=
\left\{
q\in V^*:
F(v)\geq F(w)+(q,v-w)
\text{ for all }v\in V
\right\} \\
&=
\left\{
q\in V^*:
w\in
\operatornamewithlimits{\arg\max}_{v\in V}
\{(q,v)-F(v)\}
\right\}.
\end{aligned}
\end{equation*}
If $F$ is differentiable at $w$, then
$\partial F(w)=\{F'(w)\}$.
We say that $F$ is \emph{subdifferentiable} at $w$ if
$\partial F(w)\neq\emptyset$.
The subdifferential characterization of strong convexity extends the definition of strong convexity to nondifferentiable extended real-valued functions by replacing $F'(w)$ in~\eqref{strong_convexity} with any $q \in \partial F(w)$.

\subsection{Legendre--Fenchel conjugate}

For a function
$F\colon V\to\overline{\mathbb{R}}$, its
\emph{Legendre--Fenchel conjugate}
$F^*\colon V^*\to\overline{\mathbb{R}}$ is defined as
\begin{equation}
\label{Legendre_Fenchel}
F^*(q)
=
\sup_{v\in V}\{(q,v)-F(v)\},
\quad q\in V^*.
\end{equation}
The function $F^*$ is always convex and lower semicontinuous~\cite[Theorem~12.2]{Rockafellar:1970}.

The definition~\eqref{Legendre_Fenchel} yields the
Fenchel--Young inequality
\begin{equation*}
F(v)+F^*(q) \geq (q,v), \quad v\in V,\ q\in V^*.
\end{equation*}
If $F$ is proper, convex, and lower semicontinuous, equality holds if and only if $q\in\partial F(v)$.
In this case, the subdifferentials $\partial F$ and $\partial F^*$ are inverse to each
other~\cite[Theorem~16.23]{BC:2011}:
\begin{equation}
\label{subdifferential_inverse}
q\in\partial F(v)
\quad\text{if and only if}\quad
v\in\partial F^*(q),
\quad v\in V,\ q\in V^*.
\end{equation}

If $F$ is proper, lower semicontinuous, and $\mu$-strongly convex for some $\mu>0$, then $F^*$ is
differentiable and $(F^*)'\colon V^*\to V$ is
$\mu^{-1}$-Lipschitz continuous~\cite{CP:2016a}.
In particular, if $F\colon V\to\mathbb{R}$ is $\mu$-strongly convex
and $L$-smooth, then $F^*$ is $L^{-1}$-strongly convex and
$\mu^{-1}$-smooth.
In this case, $F'\colon V\to V^*$ and
$(F^*)'\colon V^*\to V$ are inverse maps:
\begin{equation}
\label{Legendre}
\begin{aligned}
v = (F^*)' (F'(v)), \quad &v \in V, \\
q = F' ((F^*)' (q)), \quad & q \in V^*.
\end{aligned}
\end{equation}

For a convex function
$F\colon V\to\overline{\mathbb{R}}$ and a linear operator
$\Pi\colon V\to W$, the \emph{infimal postcomposition} $\Pi\triangleright F\colon W\to\overline{\mathbb{R}}$ of $F$ by $\Pi$
is defined by
\begin{equation}
\label{infimal_postcomposition}
(\Pi\triangleright F)(w)
=
\inf_{v\in V,\ \Pi v=w}F(v),
\quad w \in W.
\end{equation}
Its Legendre--Fenchel conjugate is given by (see~\cite[Proposition~13.21]{BC:2011})
\begin{equation}
\label{infimal_postcomposition_conjugate}
(\Pi\triangleright F)^*
=
F^*\circ \Pi^*.
\end{equation}

\Cref{Table:conjugate} summarizes several standard examples of Legendre--Fenchel conjugates.
The conjugate
$\chi_K^*(q)=\sup_{v\in K}(q,v)$ is called the
\emph{support function} of $K$, and
$\operatorname{LSE}_k$ denotes the \emph{log-sum-exp function};
see~\cite{BHH:2021,Murphy:2022}.
For the log-sum-exp conjugate, we use the simplex
\begin{equation}
\label{probability_simplex}
\Delta_k
=
\left\{
q\in\mathbb{R}^k:
q_i\geq0,\ 
\sum_{i=1}^k q_i=1
\right\},
\end{equation}
together with the convention $0\log 0=0$.
For further properties of the Legendre--Fenchel conjugate, we refer to~\cite{Rockafellar:1970}.

\begin{table}
\centering
\begin{tabular}{ccc}
\hline
\textbf{$F(v)$}
&
\textbf{$F^*(q)$}
&
\textbf{Notes}
\\
\hline

$\dfrac{1}{2}(Av,v)-(f,v)$
&
$\dfrac{1}{2}(q+f,A^{-1}(q+f))$
&
\makecell[l]{$A\colon V\to V^*$ SPD,\\ $f\in V^*$}
\\
\hline

$\chi_K(v)$
&
$\displaystyle\sup_{v\in K}(q,v)$
&
\makecell[l]{$K\subset V$ is nonempty\\ and convex}
\\
\hline

$\|v\|_s$
&
$\chi_{\{\|\cdot\|_{s^*}\leq1\}}(q)$
&
\makecell[l]{$V=\mathbb{R}^d$, $1\leq s\leq\infty$,\\
$1/s+1/s^*=1$}
\\
\hline

$\displaystyle
\operatorname{LSE}_k(v)
=
\log\left(\sum_{i=1}^k e^{v_i}\right)$
&
$\displaystyle
\begin{cases}
\sum_{i=1}^k q_i\log q_i,
& q\in\Delta_k,\\
\infty,
& \text{otherwise}
\end{cases}$
&
\makecell[l]{$V=\mathbb{R}^k$}
\\
\hline

\end{tabular}
\caption{Examples of Legendre--Fenchel conjugates, where $\Delta_k$ is defined in~\eqref{probability_simplex}.}
\label{Table:conjugate}
\end{table}

\begin{remark}
\label{Rem:Legendre}
Analogous inverse identities, restricted to the appropriate interiors of the effective domains, hold for the broader class of Legendre functions~\cite{Teboulle:2018}.
\end{remark}

\section{Dualization}
\label{Sec:Dualization}
In this section, we review Fenchel--Rockafellar duality and introduce the notion of dualization used throughout the paper.
We also present several motivating examples involving the augmented Lagrangian method~\cite{Hestenes:1969,Powell:1969}, the proximal point algorithm~\cite{Martinet:1970,Rockafellar:1976}, and the level-set Frank--Wolfe method~\cite{MMP:2025}.

\subsection{Fenchel--Rockafellar duality}
\label{Subsec:Fenchel_Rockafellar}

We review the key features of Fenchel--Rockafellar duality~\cite{Rockafellar:1967,Rockafellar:1970}, which establishes a duality relation between two convex optimization problems.
We revisit the structured composite convex optimization problem~\eqref{model_problem_2}, which we call the \emph{primal problem} in this section.
Here, $F\colon V\to\overline{\mathbb{R}}$ and
$G\colon W\to\overline{\mathbb{R}}$ are proper, convex, and lower semicontinuous functions, and $B\colon V\to W$ is a bounded linear operator.

Define the saddle function
$\mathcal{L}\colon V\times W^*\to\overline{\mathbb{R}}$ by
\begin{equation}
\label{saddle_function}
\mathcal{L}(v,q) = F(v) + (q,Bv) -G^*(q),
\quad v \in V, q \in W^*.
\end{equation}
By the biconjugation identity $G=G^{**}$, we obtain the weak duality relation
\begin{multline}
\label{min_max_exchange}
\inf_{v\in V}\{F(v)+G(Bv)\}
=
\inf_{v\in V}\sup_{q\in W^*}\mathcal{L}(v,q)
\\
\geq
\sup_{q\in W^*}\inf_{v\in V}\mathcal{L}(v,q)
=
\sup_{q\in W^*}
\{-F^*(-B^*q)-G^*(q)\}.
\end{multline}
We are particularly interested in the case where equality holds in~\eqref{min_max_exchange}.
A standard qualification condition is that the perturbation function
\begin{equation}
\label{Fenchel_Rockafellar_condition}
H(w)
=
\inf_{v\in V}
\{F(v)+G(Bv+w)\}
\end{equation}
be subdifferentiable at $w=0$.
Under this condition, the inequality in~\eqref{min_max_exchange} becomes an equality and the supremum on the right-hand side is attained; that is, strong duality holds~\cite[Theorem~3]{Rockafellar:1967}.
Accordingly, the Fenchel--Rockafellar \emph{dual problem} of~\eqref{model_problem_2} is defined as
\begin{equation}
\label{dual}
\min_{q\in W^*}
\{F^*(-B^*q)+G^*(q)\}.
\end{equation}
Problem~\eqref{dual} is also an instance of~\eqref{model_problem_2}.

If $u\in V$ and $p\in W^*$ are solutions of the primal and dual problems, respectively, then we have
\begin{subequations}
\label{primal_dual_relations}
\begin{align}
\label{primal_dual_relation_1}
-B^*p
&\in
\partial F(u),\\
\label{primal_dual_relation_2}
Bu
&\in
\partial G^*(p).
\end{align}
\end{subequations}
Conversely, any pair $(u,p)$ satisfying
\eqref{primal_dual_relation_1} and
\eqref{primal_dual_relation_2} solves the primal and dual problems.
We call these conditions the \emph{primal--dual relations}.
Equivalently, $(u,p)$ is a saddle point of the function
$\mathcal{L}$ defined in~\eqref{saddle_function}.

\begin{remark}
\label{Rem:Fenchel_Rockafellar_condition}
The condition~\eqref{Fenchel_Rockafellar_condition} is satisfied in many applications.
Several easily verifiable sufficient conditions can be found, e.g., in~\cite[Equation~(4.21)]{ET:1999} and~\cite[Theorems~1 and~2]{Rockafellar:1967}.
\end{remark}

Throughout the remainder of the paper, whenever Fenchel--Rockafellar duality is invoked, we assume that the corresponding strong duality condition holds, so that the associated dual formulation and primal--dual relations are valid.
We also assume that every displayed minimization subproblem admits a solution whenever an iterate is defined by an $\operatornamewithlimits{\arg\min}$ expression.

Under strict convexity of $F$, the first primal--dual relation~\eqref{primal_dual_relation_1} uniquely determines the primal solution $u$ from a dual solution $p$.
The following proposition makes this statement precise.

\begin{proposition}
\label{Prop:Fenchel_Rockafellar}
Suppose that the following hold:
\begin{enumerate}
\item[(i)] The primal problem~\eqref{model_problem_2} admits a solution.
\item[(ii)] The condition~\eqref{Fenchel_Rockafellar_condition} holds, so that the dual problem~\eqref{dual} admits a solution $p\in W^*$.
\item[(iii)] The function $F$ is strictly convex.
\end{enumerate}
If $u\in V$ satisfies~\eqref{primal_dual_relation_1}, then $u$ is the unique solution of the primal problem~\eqref{model_problem_2}.
\end{proposition}


Finally, consider the linearly constrained problem~\eqref{model_problem_3}.
Setting $G=\chi_{\{g\}}$ in~\eqref{model_problem_2} gives~\eqref{model_problem_3}.
Its Fenchel--Rockafellar dual problem is written as
\begin{equation}
\label{PPA_dual}
\min_{q \in W^*} \left\{ F^* (-B^* q) + (q,g) \right\}.
\end{equation}
The corresponding primal--dual relations are
\begin{equation}
\label{PPA_primal_dual_relations}
-B^*p\in\partial F(u),
\quad Bu=g.
\end{equation}
Problem~\eqref{PPA_dual} is an instance of~\eqref{model_problem_2}, which in turn is a special case of~\eqref{model_problem_1}.
Thus, Fenchel--Rockafellar duality transforms~\eqref{model_problem_3} into a problem of the form~\eqref{model_problem_1}.
Consequently, through specialization and Fenchel--Rockafellar duality, the three model problems~\eqref{model_problem_1}--\eqref{model_problem_3} can all be treated within the framework of~\eqref{model_problem_1}.

\subsection{Dualization of convex optimization algorithms}
We now give the formal definition of dualization, motivated by \cref{Prop:Fenchel_Rockafellar}.
Since an iterative method generates a sequence of iterates, we refer to the method through the sequence it produces when no confusion can arise.

\begin{definition}
\label{Def:dualization}
An iterative method $\{ u^{(n)} \} \subset V$ for solving the primal problem~\eqref{model_problem_2} is said to be a dualization of an iterative method $\{ p^{(n)} \}\subset W^*$ for solving the dual problem~\eqref{dual} if, for any initial configuration of the dual algorithm, there exists an initial configuration of the primal algorithm such that
\begin{equation}
\label{Def1:dualization}
- B^* p^{(n)} \in \partial F (u^{(n)}), \quad n \geq 1.
\end{equation}
\end{definition}

The definition is motivated by the first primal--dual relation~\eqref{primal_dual_relation_1}.
It formalizes a mechanism implicit in classical primal--dual algorithmic correspondences.
In the examples below, corresponding local subproblems are related through Fenchel--Rockafellar duality, and the resulting iterates are shown to satisfy~\eqref{Def1:dualization}, which permits convergence and structural properties to be transferred between the two formulations.

\begin{remark}
\label{Rem:dualization}
By duality, one may instead transform a method for the primal problem~\eqref{model_problem_2} into a method for the dual problem~\eqref{dual} by requiring
\begin{equation*}
B u^{(n)} \in \partial G^* (p^{(n)}), \quad n \geq 1,
\end{equation*}
as suggested by the second primal--dual relation~\eqref{primal_dual_relation_2}.
Applying dualization twice need not recover the original method, since the two primal--dual relations need not hold simultaneously at general iterates \((u^{(n)},p^{(n)})\).
\end{remark}

Under strong convexity of $F$, the following theorem shows that convergence of the dual sequence implies convergence of its dualization; in particular, a linear rate is preserved up to a multiplicative constant.

\begin{theorem}
\label{Thm:convergence_transfer}
Suppose that the following hold:
\begin{enumerate}
\item[(i)] The primal problem~\eqref{model_problem_2} admits a solution $u \in V$.
\item[(ii)] The condition~\eqref{Fenchel_Rockafellar_condition} holds, ensuring that the dual problem~\eqref{dual} admits a solution $p \in W^*$.
\item[(iii)] The function $F$ in~\eqref{model_problem_2} is $\mu$-strongly convex for some $\mu > 0$.
\end{enumerate}
If an iterative method $\{ u^{(n)} \}$ for solving the primal problem~\eqref{model_problem_2} is a dualization of an iterative method $\{ p^{(n)} \}$ for solving the dual problem~\eqref{dual}, then we have
\begin{equation*}
\| u^{(n)} - u \| \leq \mu^{-1} \| B^* (p^{(n)} - p) \|, \quad n \geq 1.
\end{equation*}
In particular, if $p^{(n)} \to p$ in $W^*$, then $u^{(n)} \to u$ in $V$.
\end{theorem}
\begin{proof}
Since $F$ is $\mu$-strongly convex, $(F^*)'$ is $\mu^{-1}$-Lipschitz continuous~\cite{CP:2016a}.
Moreover, by~\eqref{subdifferential_inverse}, we see that~\eqref{Def1:dualization} implies
\begin{equation*}
u^{(n)} = (F^*)' (-B^* p^{(n)}), \quad n \geq 1.
\end{equation*}
Similarly,~\eqref{primal_dual_relation_1} implies
\begin{equation*}
u = (F^*)' (-B^* p).
\end{equation*}
It follows that
\begin{equation*}
\| u^{(n)} - u \|
= \| (F^*)' (-B^* p^{(n)}) - (F^*)' (-B^* p) \|
\leq \mu^{-1} \| B^* (p^{(n)} - p) \|,
\end{equation*}
which completes the proof.
\end{proof}

\begin{remark}
\label{Rem:convergence_transfer}
The conditions in \cref{Thm:convergence_transfer} that ensure convergence of the primal sequence \(u^{(n)}\) can be weakened.
If \(F\) is uniformly convex, then the same monotonicity argument gives a convergence estimate in terms of the modulus of uniform convexity; in the special case of power-type uniform convexity, this yields a power-type bound in terms of \( \|B^*(p^{(n)}-p)\| \).
If \(F\) is strictly convex, the same conclusion follows if \(\{u^{(n)}\}\) is relatively compact, since every limit point must satisfy \(-B^* p\in\partial F(v)\), which has at most one solution.
We omit the details here and keep the strongly convex case for clarity of exposition.
\end{remark}

\subsection{Augmented Lagrangian method and its dualization}

The augmented Lagrangian method~\cite{Hestenes:1969,Powell:1969} is a fundamental algorithm for solving the linearly constrained problem~\eqref{model_problem_3}.
An equivalent saddle point formulation is
\begin{equation}
\label{ALM_saddle}
\min_{v \in V} \max_{\lambda \in W^*}
\left\{ F(v) + (\lambda, Bv-g) \right\}.
\end{equation}
Given a penalty parameter $\beta>0$, define the augmented Lagrangian by
\begin{equation}
\label{ALM_augmented_Lagrangian}
\mathcal{L}_{\beta}(v,\lambda)
=
F(v)+(\lambda,Bv-g)+\frac{\beta}{2}\|Bv-g\|^2.
\end{equation}
Since the supremum with respect to $\lambda$ enforces $Bv=g$, problem~\eqref{ALM_saddle} is equivalent to the augmented Lagrangian formulation
\begin{equation*}
\min_{v \in V} \max_{\lambda \in W^*}
\mathcal{L}_{\beta}(v,\lambda).
\end{equation*}
The augmented Lagrangian method alternates minimization with respect to $v$ and a gradient ascent step with respect to $\lambda$ as follows:
\begin{subequations}
\label{Alg:ALM}
\begin{align}
    \label{Alg1:ALM}
    u^{(n+1)}
    &\in \operatornamewithlimits{\arg\min}_{v \in V}
    \left\{
        F(v) + ( \lambda^{(n)}, Bv-g)
        + \frac{\beta}{2}\|Bv-g\|^2
    \right\}, \\
    \label{Alg2:ALM}
    \lambda^{(n+1)}
    &= \lambda^{(n)} + \beta(Bu^{(n+1)}-g).
\end{align}
\end{subequations}

We next relate the augmented Lagrangian method to the proximal point algorithm~\cite{Martinet:1970,Rockafellar:1976}.
As shown above, the Fenchel--Rockafellar dual of~\eqref{model_problem_3} is~\eqref{PPA_dual}, and the corresponding primal--dual relations are given in~\eqref{PPA_primal_dual_relations}.
Given the same parameter $\beta>0$, the proximal point algorithm for solving~\eqref{PPA_dual} is written as
\begin{equation}
\label{Alg:PPA_dual}
    p^{(n+1)}
    = \operatornamewithlimits{\arg\min}_{q \in W^*}
    \left\{
        F^*(-B^* q) + (q,g)
        + \frac{1}{2\beta} \|q-p^{(n)}\|^2
    \right\}.
\end{equation}
We now recover the well-known equivalence~\cite{EB:1992,Setzer:2011,WT:2010} by proving that the augmented Lagrangian method~\eqref{Alg:ALM} is a dualization of the proximal point algorithm~\eqref{Alg:PPA_dual}.
Suppose that $\lambda^{(n)} = p^{(n)}$ holds for some $n \geq 0$.
The $p^{(n+1)}$-subproblem~\eqref{Alg:PPA_dual} is an instance of~\eqref{model_problem_2} with
\begin{equation*}
\resizebox{\textwidth}{!}{$\displaystyle
V \leftarrow W^*,\quad
W \leftarrow V^*,\quad
F(q) \leftarrow \frac{1}{2\beta}\|q-p^{(n)}\|^2 + (q,g),\quad
G(w) \leftarrow F^*(w),\quad
B \leftarrow -B^*.
$}
\end{equation*}
Invoking Fenchel--Rockafellar duality, we obtain the following dual formulation of the $p^{(n+1)}$-subproblem:
\begin{equation}
\label{PPA_subproblem_dual}
    \min_{v \in V}
    \left\{
        F(v) + ( p^{(n)}, Bv-g)
        + \frac{\beta}{2}\|Bv-g\|^2
    \right\}.
\end{equation}
Under the assumption $\lambda^{(n)}=p^{(n)}$, the problem~\eqref{PPA_subproblem_dual} agrees with~\eqref{Alg1:ALM}.
Thus, $u^{(n+1)}$ is chosen as a solution of~\eqref{PPA_subproblem_dual}.
The first primal--dual relation for this subproblem reads as
\begin{equation}
    \label{PPA_subproblem_relation1}
    Bu^{(n+1)} = g + \frac{1}{\beta} (p^{(n+1)}-p^{(n)} ).
\end{equation}
By~\eqref{PPA_subproblem_relation1}, we have
\begin{equation*}
    p^{(n+1)} = p^{(n)} + \beta(Bu^{(n+1)}-g).
\end{equation*}
Comparing this identity with the multiplier update~\eqref{Alg2:ALM}, we obtain $\lambda^{(n+1)}=p^{(n+1)}$.
Therefore, if $\lambda^{(0)}=p^{(0)}$, then by mathematical induction,
\begin{equation*}
    \lambda^{(n)}=p^{(n)}, \quad
    -B^* p^{(n)} \in \partial F(u^{(n)}),
    \quad n \geq 1.
\end{equation*}
The latter relation is precisely the defining relation~\eqref{Def1:dualization}.
Hence, the augmented Lagrangian method for~\eqref{model_problem_3} is a dualization, in the sense of \cref{Def:dualization}, of the proximal point algorithm applied to the dual problem~\eqref{PPA_dual}.

\begin{remark}
\label{Rem:PPA}
Since the split Bregman method, which is commonly used in applications such as total variation denoising and compressed sensing~\cite{COS:2010,GO:2009}, is known to generate the same primal sequence as the augmented Lagrangian method~\cite{WT:2010}, we also conclude that the split Bregman method is a dualization of the proximal point algorithm.
\end{remark}

\subsection{Dualized level-set Frank--Wolfe method}

As an example from recent literature, we interpret the dualized level-set Frank--Wolfe method of~\cite{MMP:2025}, which is based on the classical Frank--Wolfe algorithm~\cite{FW:1956}, within our dualization framework.
Given convex sets $K\subset V$ and $Q\subset W^*$, consider the constrained saddle point problem
\begin{equation}
\label{MMP_saddle}
\min_{v\in K}\max_{\lambda\in Q}
\left\{F(v)+(\lambda,Bv-g)\right\}.
\end{equation}
This formulation extends the saddle point problem~\eqref{ALM_saddle}.
Two choices of $Q$ are given as follows:
\begin{enumerate}
\item
If $W=\mathbb{R}^m$ and
$Q=\{\lambda\in\mathbb{R}^m:\lambda\leq0\}$, then the corresponding primal formulation represents the componentwise inequality constraint $Bv-g\geq0$.
\item
Bounded choices of $Q$ yield penalty terms; for example,
$Q=\{\lambda\in\mathbb{R}^m:\|\lambda\|_{\infty}\leq1\}$ yields the penalty $\|Bv-g\|_1$.
\end{enumerate}
In~\cite{MMP:2025}, the following equivalent formulation of~\eqref{MMP_saddle} is considered:
\begin{equation}
\label{MMP_primal}
\min_{v\in V}
\left\{
F(v)+\chi_K(v)+\chi_Q^*(Bv-g)
\right\},
\end{equation}
where $F\colon V\to\mathbb{R}$ is convex with Lipschitz continuous derivative, $K\subset V$ is a nonempty compact convex set, $Q\subset W^*$ is a nonempty closed convex set, $B\colon V\to W$ is a bounded linear operator, and $g\in W$.
For additional structural, qualification, and implementability assumptions, we refer to~\cite{MMP:2025}.

Since~\eqref{MMP_primal} is a special case of~\eqref{model_problem_2}, we first derive its dual problem.
In~\eqref{model_problem_2}, set
\begin{equation*}
\resizebox{\textwidth}{!}{$\displaystyle
W \leftarrow V\times W,\quad
B \leftarrow \begin{bmatrix} -I \\ 0 \end{bmatrix},\quad
G(z,w) \leftarrow \chi_{-K}(z)+\chi_Q^*(w-Bz-g),
\quad (z,w)\in V\times W,
$}
\end{equation*}
with $F$ unchanged, where the $B$ appearing in the definition of $G$ denotes the operator in~\eqref{MMP_primal}.
Let
\begin{equation*}
K_B=\{(v,Bv):v\in K\}\subset V\times W.
\end{equation*}
For $(q,\mu)\in V^*\times W^*$, we have
\begin{equation*}
\begin{bmatrix} -I \\ 0 \end{bmatrix}^*(q,\mu)=-q,\quad
G^*(q,\mu)
=
\chi_Q(\mu)+(\mu,g)
+\chi_{-K_B}^*(q,\mu).
\end{equation*}
Hence, Fenchel--Rockafellar duality yields the dual problem
\begin{equation}
\label{MMP_dual}
\min_{(q,\mu)\in V^*\times W^*}
\left\{
F^*(q)
+\chi_Q(\mu)
+(\mu,g)
+\chi_{-K_B}^*(q,\mu)
\right\}.
\end{equation}
Let $u$ be a primal solution of~\eqref{MMP_primal}, and let
$(p,\lambda)\in V^*\times W^*$ be a dual solution of~\eqref{MMP_dual}.
The corresponding first primal--dual relation is
\begin{equation}
\label{MMP_primal_dual_relation}
p\in\partial F(u).
\end{equation}

We next consider a projection method for the dual problem~\eqref{MMP_dual}.
At each iteration, the method constructs a closed convex approximation
$\mathcal{Q}^{(n)}\subset V^*\times W^*$ of a level set of the dual objective from the information generated in the previous iterations, and computes the next dual iterate by a Bregman projection onto $\mathcal{Q}^{(n)}$.
This method is called the extended level-set projection method in~\cite{MMP:2025}.
Suppose that
$(p^{(n)},\lambda^{(n)})\in V^*\times W^*$
and $u^{(n)}\in V$ satisfy
\begin{equation*}
p^{(n)}\in\partial F(u^{(n)}).
\end{equation*}
The corresponding Bregman projection subproblem is
\begin{equation}
\label{MMP_dual_update}
\begin{split}
(p^{(n+1)},\lambda^{(n+1)})
\in
\operatornamewithlimits{\arg\min}_{(q,\mu)\in\mathcal{Q}^{(n)}}
\left\{
F^*(q)-F^*(p^{(n)})
-(q-p^{(n)},u^{(n)})
+\frac{1}{2}\|\mu-\lambda^{(n)}\|^2
\right\}.
\end{split}
\end{equation}
For details on the construction of $\mathcal{Q}^{(n)}$, we refer to~\cite{MMP:2025}.

We next dualize the subproblem~\eqref{MMP_dual_update}.
Ignoring constants independent of $(q,\mu)$, it is written as
\begin{equation}
\label{MMP_dual_update_rewritten}
(p^{(n+1)},\lambda^{(n+1)})
\in
\operatornamewithlimits{\arg\min}_{(q,\mu)\in V^*\times W^*}
\left\{
F^*(q)-(q,u^{(n)})
+\frac{1}{2}\|\mu-\lambda^{(n)}\|^2
+\chi_{\mathcal{Q}^{(n)}}(q,\mu)
\right\}.
\end{equation}
This is an instance of~\eqref{model_problem_2}, where the linear operator is the projection $(q,\mu)\mapsto q$.
Applying Fenchel--Rockafellar duality to~\eqref{MMP_dual_update_rewritten}, we obtain the following primal-side subproblem:
\begin{equation}
\label{MMP_primal_update}
u^{(n+1)}
\in
\operatornamewithlimits{\arg\min}_{v\in V}
\left\{
F(v)
+
\sup_{(q,\mu)\in\mathcal{Q}^{(n)}}
\left\{
(q,u^{(n)}-v)
-\frac{1}{2}\|\mu-\lambda^{(n)}\|^2
\right\}
\right\}.
\end{equation}
This primal-side subproblem underlies the Frank--Wolfe-type update in~\cite{MMP:2025}, which uses finite information generated by linear minimization over $K$.
The primal--dual relations for the pair of subproblems~\eqref{MMP_dual_update_rewritten} and~\eqref{MMP_primal_update} give
\begin{equation*}
u^{(n+1)}\in\partial F^*(p^{(n+1)}),
\end{equation*}
or equivalently, by~\eqref{subdifferential_inverse},
\begin{equation*}
p^{(n+1)}\in\partial F(u^{(n+1)}),
\end{equation*}
which has the same form as~\eqref{MMP_primal_dual_relation}.
Hence, the primal-side method obtained from~\eqref{MMP_primal_update} is a dualization of the projection method~\eqref{MMP_dual_update} in view of \cref{Def:dualization}.

\section{Subspace correction methods}
\label{Sec:MSC}
Subspace correction methods~\cite{Xu:1992,XZ:2002} provide a general framework for a broad class of iterative methods, ranging from classical Jacobi and Gauss--Seidel methods to domain decomposition and multigrid methods~\cite{Bramble:1993,BPX:1990,TW:2005}.
In this section, we first review their formulation for linear systems and their equivalence to block methods for an expanded system.
We then extend the same constructions to convex optimization, which is the setting used in the remainder of this paper.

\subsection{Space decomposition}
Our fundamental assumption is that $V$ admits a decomposition
\begin{equation}
\label{MSC_space_decomposition}
V=\sum_{j=1}^J \Pi_j V_j,
\end{equation}
where each $V_j$, $1\leq j\leq J$, is a real Hilbert space, referred to as a \emph{local space}, and $\Pi_j\colon V_j\to V$ is a bounded linear operator.
If $V_j$ is a subspace of $V$ and $\Pi_j$ is the natural embedding, then this reduces to the standard subspace decomposition considered in much of the literature~\cite{Xu:1992,XZ:2002}.

For later use, we define the product space $\undertilde{V} = \prod_{j=1}^J V_j$ and the summation operator $\Pi \colon \undertilde{V} \to V$ by
\begin{equation*}
\Pi\undertilde{v}=\sum_{j=1}^J\Pi_jv_j,
\quad \undertilde{v} = (v_1, \dots, v_J) \in \undertilde{V}.
\end{equation*}
It follows from the space decomposition~\eqref{MSC_space_decomposition} that $\Pi$ is surjective.
Different choices of the local spaces and transfer operators lead to block, domain decomposition, and multigrid methods; see~\cite{TX:2002,TW:2005,Xu:1992}.

\subsection{Linear systems}
\label{Subsec:MSC_linear}
We first consider the linear system
\begin{equation}
\label{MSC_linear_system}
Au=f,
\end{equation}
where $A\colon V\to V^*$ is SPD and $f\in V^*$.
We denote the unique solution of~\eqref{MSC_linear_system} by $u\in V$.
We use the energy inner product and norm
\begin{equation*}
(v,w)_A=(Av,w),
\quad
\|v\|_A=(Av,v)^{1/2}.
\end{equation*}

For each local space $V_j$, we define the local operator
\begin{equation*}
A_j=\Pi_j^*A\Pi_j\colon V_j\to V_j^*.
\end{equation*}
For simplicity, we assume that each $A_j$ is SPD; the semidefinite case is treated in~\cite{LWXZ:2008,PX:2025}.
We choose a boundedly invertible local solver $R_j\colon V_j^*\to V_j$ that approximates $A_j^{-1}$.
Given a current approximation $u^{\mathrm{old}} \in V$, the $j$th local correction is given by
\begin{equation}
\label{MSC_linear_local}
w_j=R_j\Pi_j^*(f-Au^{\mathrm{old}}),
\end{equation}
and the corresponding local update is $u^{\mathrm{new}}=u^{\mathrm{old}}+\Pi_jw_j$.
When $R_j=A_j^{-1}$, the correction~\eqref{MSC_linear_local} is the exact solution of the local residual equation
\begin{equation*}
(A_jw_j,z_j)=(\Pi_j^*(f-Au^{\mathrm{old}}),z_j),
\quad z_j\in V_j.
\end{equation*}

The \emph{parallel subspace correction method}~(PSC) computes all local corrections from the same iterate and combines them simultaneously, as presented in \cref{Alg:MSC_linear_PSC}.
PSC is also known as the additive Schwarz method in domain decomposition and the Bramble--Pasciak--Xu preconditioner in multigrid methods~\cite{BPX:1990,TW:2005,Xu:1992}.

\begin{algorithm}
\caption{Parallel subspace correction method~(PSC) for~\eqref{MSC_linear_system}}
\label{Alg:MSC_linear_PSC}
\begin{algorithmic}[]
\STATE Given $\tau>0$ and choose $u^{(0)}\in V$.
\FOR{$n=0,1,2,\dots$}
    \FOR{$j=1,2,\dots,J$ \textbf{in parallel}}
        \STATE \vspace{-0.5cm}
        \begin{equation*}
        w_j^{(n+1)}
        =R_j\Pi_j^*(f-Au^{(n)})
        \end{equation*}
        \vspace{-0.3cm}
    \ENDFOR
    \STATE \vspace{-0.5cm}
    \begin{equation*}
    u^{(n+1)}
    =u^{(n)}+\tau\sum_{j=1}^J\Pi_jw_j^{(n+1)}
    \end{equation*}
    \vspace{-0.3cm}
\ENDFOR
\end{algorithmic}
\end{algorithm}

The \emph{successive subspace correction method}~(SSC) updates one local space at a time, as presented in \cref{Alg:MSC_linear_SSC}.
SSC corresponds to the multiplicative Schwarz method in domain decomposition and the basic \textbackslash-cycle in multigrid methods~\cite{BTW:2003,Bramble:1993,TW:2005}.

\begin{algorithm}
\caption{Successive subspace correction method~(SSC) for~\eqref{MSC_linear_system}}
\label{Alg:MSC_linear_SSC}
\begin{algorithmic}[]
\STATE Choose $u^{(0)}\in V$.
\FOR{$n=0,1,2,\dots$}
    \FOR{$j=1,2,\dots,J$}
        \STATE \vspace{-0.5cm}
        \begin{align*}
        w_j^{(n+1)}
        &=R_j\Pi_j^* (f-Au^{(n+\frac{j-1}{J})} ),\\
        u^{(n+\frac{j}{J})}
        &=u^{(n+\frac{j-1}{J})}+\Pi_jw_j^{(n+1)}.
        \end{align*}
        \vspace{-0.3cm}
    \ENDFOR
\ENDFOR
\end{algorithmic}
\end{algorithm}

We next present sharp identities for PSC and SSC that are useful for convergence rate analysis~\cite{Brenner:2013,PX:2025,XZ:2002}.

\begin{theorem}
\label{Thm:PSC_convergence}
Suppose that each local solver $R_j\colon V_j^*\to V_j$ is SPD, and define
\begin{equation*}
B_{\mathrm{PSC}}
=
\sum_{j=1}^J\Pi_jR_j\Pi_j^*
\colon V^*\to V.
\end{equation*}
Then $B_{\mathrm{PSC}}$ is SPD and
\begin{equation*}
(B_{\mathrm{PSC}}^{-1}v,v)
=
\inf_{\sum_{j=1}^J\Pi_jv_j=v}
\sum_{j=1}^J(R_j^{-1}v_j,v_j),
\quad v\in V.
\end{equation*}
\end{theorem}

\begin{theorem}[Xu--Zikatanov identity]
\label{Thm:SSC_convergence}
For $1\leq j\leq J$, suppose that $R_j\colon V_j^*\to V_j$ is boundedly invertible and that
\begin{equation*}
\overline{R}_j
=
R_j+R_j^*-R_j^*A_jR_j
\colon V_j^*\to V_j
\end{equation*}
is SPD.
Define $B_{\mathrm{SSC}}\colon V^*\to V$ by
\begin{equation*}
I-B_{\mathrm{SSC}}A
=
(I-\Pi_JR_J\Pi_J^*A)\cdots(I-\Pi_1R_1\Pi_1^*A).
\end{equation*}
Then SSC converges and satisfies
\begin{equation*}
\|I-B_{\mathrm{SSC}}A\|_A^2
=
1-\frac{1}{1+c_0}
<1,
\end{equation*}
where
\begin{equation*}
c_0
=
\sup_{\|v\|_A=1}
\inf_{\sum_{j=1}^J\Pi_jv_j=v}
\sum_{j=1}^J
\|R_j^*w_j\|_{\overline{R}_j^{-1}}^2,
\quad w_i
=
\Pi_i^*A\sum_{j=i}^J\Pi_jv_j-R_i^{-1}v_i,
\quad 1\leq i\leq J.
\end{equation*}
\end{theorem}

\begin{remark}
\label{Rem:aux_space_lemma}
The infimal postcomposition identity~\eqref{infimal_postcomposition_conjugate}
is closely related to \cref{Thm:PSC_convergence}.
Indeed, let
\begin{equation}
\label{R_undertilde}
\undertilde{R}=\operatorname{diag}(R_1,\dots,R_J)\colon\undertilde{V}^*\to\undertilde{V}.
\end{equation}
For
\begin{equation*}
F(\undertilde{v})
=\frac{1}{2}
(\undertilde{R}^{-1}\undertilde{v},\undertilde{v}),
\quad \undertilde{v}\in\undertilde{V},
\end{equation*}
the identity~\eqref{infimal_postcomposition_conjugate} and the definition of
$B_{\mathrm{PSC}}$ give
\begin{equation*}
(\Pi\triangleright F)^*(q)
=F^*(\Pi^*q)
=\frac{1}{2}(\undertilde{R}\Pi^*q,\Pi^*q)
=\frac{1}{2}(B_{\mathrm{PSC}}q,q),
\quad q\in V^*.
\end{equation*}
Since $\Pi$ is surjective and $\undertilde{R}$ is SPD,
$\Pi\triangleright F$ is continuous and convex.
Taking the Legendre--Fenchel conjugate once more yields
\begin{equation*}
(\Pi\triangleright F)(v)
=\frac{1}{2}(B_{\mathrm{PSC}}^{-1}v,v),
\quad v\in V,
\end{equation*}
which is equivalent to the identity in \cref{Thm:PSC_convergence}.
\end{remark}

There is a relationship between SSC and the classical alternating projection method known as the von Neumann algorithm~\cite{Neumann:1950}.
Let $M_j$, $1\leq j\leq J$, be closed subspaces of $V$, assume that $\sum_{j=1}^J M_j^\perp$ is closed, and set $M=\bigcap_{j=1}^J M_j$.
Given $f\in V$, the von Neumann alternating projection algorithm for computing $\operatorname{proj}_M f$ is presented in \cref{Alg:MSC_von_Neumann}.

\begin{algorithm}
\caption{von Neumann alternating projection algorithm}
\label{Alg:MSC_von_Neumann}
\begin{algorithmic}[]
\STATE Set $u^{(0)}=f$.
\FOR{$n=0,1,2,\dots$}
    \FOR{$j=1,2,\dots,J$}
        \STATE \vspace{-0.5cm}
        \begin{equation*}
        u^{(n+\frac{j}{J})}
        =
        \operatorname{proj}_{M_j}u^{(n+\frac{j-1}{J})}
        \end{equation*}
        \vspace{-0.3cm}
    \ENDFOR
\ENDFOR
\end{algorithmic}
\end{algorithm}

Since $M^\perp=\sum_{j=1}^J M_j^\perp$, the projection $p=\operatorname{proj}_{M^\perp}f$ is characterized by the variational problem
\begin{equation}
\label{von_Neumann_variational}
(p,s)=(f,s),
\quad s\in M^\perp.
\end{equation}
Applying SSC to~\eqref{von_Neumann_variational} with the local spaces $M_j^\perp$, the exact local solvers induced by their inner products, and $p^{(0)}=0$, we obtain
\begin{align*}
r_j^{(n+1)}
&=\operatorname{proj}_{M_j^\perp}
(f-p^{(n+\frac{j-1}{J})}),\\
p^{(n+\frac{j}{J})}
&=p^{(n+\frac{j-1}{J})}+r_j^{(n+1)}.
\end{align*}
If we define $u^{(n+\frac{j}{J})}=f-p^{(n+\frac{j}{J})}$, then
\begin{equation*}
u^{(n+\frac{j}{J})}
=
u^{(n+\frac{j-1}{J})}
-
\operatorname{proj}_{M_j^\perp}
u^{(n+\frac{j-1}{J})}
=
\operatorname{proj}_{M_j}
u^{(n+\frac{j-1}{J})}.
\end{equation*}
Hence, under the relation $u^{(n)}=f-p^{(n)}$, SSC is exactly the von Neumann algorithm in~\cref{Alg:MSC_von_Neumann}; this equivalence was established in~\cite{XZ:2002}.

\subsection{Block methods for the expanded system}
Based on the space decomposition~\eqref{MSC_space_decomposition}, we consider the following expanded linear system
\begin{equation}
\label{MSC_linear_expanded}
\undertilde{A}\undertilde{u}=\undertilde{f},
\quad \text{ where }
\undertilde{A}=\Pi^*A\Pi,\
\undertilde{f}=\Pi^*f.
\end{equation}
Problem~\eqref{MSC_linear_expanded} is the optimality system for a quadratic instance of~\eqref{model_problem_1} on the product space $\undertilde{V}$.
The operator $\undertilde{A}\colon\undertilde{V}\to\undertilde{V}^*$ has the block form
\begin{equation*}
\undertilde{A}=[A_{ij}]_{i,j=1}^J,
\quad
A_{ij}=\Pi_i^*A\Pi_j,
\end{equation*}
We also write $\undertilde{f}=(f_1,\dots,f_J)$, where $f_j=\Pi_j^*f$.
Then $A_{jj}=A_j$.
Since $\Pi$ is surjective, a vector $\undertilde{u}\in\undertilde{V}$ solves~\eqref{MSC_linear_expanded} if and only if $u=\Pi\undertilde{u}$ solves~\eqref{MSC_linear_system}~\cite[Proposition~3.2]{PX:2025}.
The expanded system may be singular even though $A$ is SPD, because $\undertilde{A}$ has the null space $\ker\Pi$.

The expanded Jacobi and expanded Gauss--Seidel iterations for~\eqref{MSC_linear_expanded} are presented in \cref{Alg:MSC_linear_block_Jacobi,Alg:MSC_linear_block_GS}, respectively.

\begin{algorithm}
\caption{Expanded Jacobi method for~\eqref{MSC_linear_expanded}}
\label{Alg:MSC_linear_block_Jacobi}
\begin{algorithmic}[]
\STATE Given $\tau>0$ and choose $\undertilde{u}^{(0)}\in\undertilde{V}$.
\FOR{$n=0,1,2,\dots$}
    \FOR{$j=1,2,\dots,J$ \textbf{in parallel}}
        \STATE \vspace{-0.5cm}
        \begin{equation*}
        u_j^{(n+1)}
        =u_j^{(n)}
        +\tau R_j
        \left(
        f_j-\sum_{i=1}^J A_{ji}u_i^{(n)}
        \right)
        \end{equation*}
        \vspace{-0.3cm}
    \ENDFOR
\ENDFOR
\end{algorithmic}
\end{algorithm}

\begin{algorithm}
\caption{Expanded Gauss--Seidel method for~\eqref{MSC_linear_expanded}}
\label{Alg:MSC_linear_block_GS}
\begin{algorithmic}[]
\STATE Choose $\undertilde{u}^{(0)}\in\undertilde{V}$.
\FOR{$n=0,1,2,\dots$}
    \FOR{$j=1,2,\dots,J$}
        \STATE \vspace{-0.5cm}
        \begin{equation*}
        u_j^{(n+1)}
        =u_j^{(n)}
        +R_j\left(
        f_j
        -\sum_{i=1}^{j-1}A_{ji}u_i^{(n+1)}
        -\sum_{i=j}^J A_{ji}u_i^{(n)}
        \right)
        \end{equation*}
        \vspace{-0.3cm}
    \ENDFOR
\ENDFOR
\end{algorithmic}
\end{algorithm}

The linear PSC and SSC iterations in \cref{Alg:MSC_linear_PSC,Alg:MSC_linear_SSC} are equivalent to the expanded Jacobi and expanded Gauss--Seidel iterations in \cref{Alg:MSC_linear_block_Jacobi,Alg:MSC_linear_block_GS}, respectively.
More precisely, if $u^{(0)}=\Pi\undertilde{u}^{(0)}$, then the corresponding iterates satisfy
\begin{equation*}
u^{(n)}=\Pi\undertilde{u}^{(n)},
\quad n\geq1.
\end{equation*}
This equivalence is a direct consequence of the auxiliary space formulation; see~\cite{Chen:2011} and~\cite[Section~6]{PX:2025} for details.
Thus, subspace correction methods for a possibly overlapping or auxiliary decomposition can be viewed as elementary block methods on a larger product space.

\subsection{Extension to convex optimization}
Subspace correction methods were originally developed for linear problems, and several approaches have extended them to nonlinear problems.
The Newton--Krylov--Schwarz method uses a Newton-type outer iteration and applies subspace correction to the resulting linearized problems; see~\cite{CGKMY:1998,Keyes:1995,KK:2004}.
Related applications to variational inequalities can be found in~\cite{BDS:2018,BSW:2022,Lee:2013}.
Another approach is nonlinear preconditioning~\cite{CK:2002,DGKKM:2016,LKK:2018}.
Overlapping Schwarz constructions were studied in~\cite{CKM:2002,HC:2005,HC:2007}, with a recent development in~\cite{LGYK:2024}.
In the following, we consider a direct extension of subspace correction methods to convex optimization.

We consider the convex optimization problem~\eqref{model_problem_1},
where $E\colon V\to\overline{\mathbb{R}}$ is proper, convex, lower semicontinuous, and coercive.
Then~\eqref{model_problem_1} admits a solution $u\in V$~\cite[Proposition~11.14]{BC:2011}.

Given a current approximation $u^{\mathrm{old}}\in V$, an exact correction in the $j$th local space is obtained, whenever the following problem admits a solution, by solving
\begin{equation}
\label{MSC_local_exact}
\min_{w_j\in V_j}E(u^{\mathrm{old}}+\Pi_jw_j).
\end{equation}
To allow inexact and computationally convenient local solvers, we instead consider
\begin{equation}
\label{MSC_local}
\min_{w_j\in V_j}E_j(w_j;u^{\mathrm{old}}),
\end{equation}
where $E_j(w_j;u^{\mathrm{old}})$ is a suitable approximation of $E(u^{\mathrm{old}}+\Pi_jw_j)$.
We assume that~\eqref{MSC_local} has a unique solution for every $u^{\mathrm{old}}\in V$.
Rigorous approximability conditions and examples of inexact local models can be found in~\cite[Assumption~4.3 and Section~6.4]{Park:2020} and~\cite[Assumption~A.1]{LP:2025b}.

The parallel and successive extensions of the linear methods are presented in \cref{Alg:PSC,Alg:SSC}, respectively.

\begin{algorithm}
\caption{Parallel subspace correction method~(PSC) for~\eqref{model_problem_1}}
\label{Alg:PSC}
\begin{algorithmic}[]
\STATE Given $\tau>0$ and choose $u^{(0)}\in V$.
\FOR{$n=0,1,2,\dots$}
    \FOR{$j=1,2,\dots,J$ \textbf{in parallel}}
        \STATE \vspace{-0.5cm}\begin{equation*}
        w_j^{(n+1)}\in
        \operatornamewithlimits{\arg\min}_{w_j\in V_j}
        E_j(w_j;u^{(n)})
        \end{equation*}
        \vspace{-0.3cm}
    \ENDFOR
    \STATE \vspace{-0.5cm}\begin{equation*}
    u^{(n+1)}=u^{(n)}+\tau\sum_{j=1}^J\Pi_jw_j^{(n+1)}
    \end{equation*}
    \vspace{-0.3cm}
\ENDFOR
\end{algorithmic}
\end{algorithm}

\begin{algorithm}
\caption{Successive subspace correction method~(SSC) for~\eqref{model_problem_1}}
\label{Alg:SSC}
\begin{algorithmic}[]
\STATE Choose $u^{(0)}\in V$.
\FOR{$n=0,1,2,\dots$}
    \FOR{$j=1,2,\dots,J$}
        \STATE \vspace{-0.5cm}\begin{align*}
        w_j^{(n+1)}
        &\in\operatornamewithlimits{\arg\min}_{w_j\in V_j}
        E_j (w_j;u^{(n+\frac{j-1}{J})} ),\\
        u^{(n+\frac{j}{J})}
        &=u^{(n+\frac{j-1}{J})}+\Pi_jw_j^{(n+1)}.
        \end{align*}
        \vspace{-0.3cm}
    \ENDFOR
\ENDFOR
\end{algorithmic}
\end{algorithm}

\begin{remark}
\label{Rem:linear}
Consider the quadratic energy
\begin{equation*}
E(v)=\frac{1}{2}(Av,v)-(f,v),
\quad v\in V,
\end{equation*}
where $A\colon V\to V^*$ is SPD and $f\in V^*$.
Then~\eqref{model_problem_1} is equivalent to the linear system~\eqref{MSC_linear_system}, which admits a unique solution $u\in V$.
For this particular convex local model, suppose that $R_j^{-1}$ induces an SPD bilinear form and define
\begin{equation}
\label{MSC_quadratic_local_model}
E_j(w_j;v)
=E(v)+(E'(v),\Pi_jw_j)
+\frac{1}{2}(R_j^{-1}w_j,w_j).
\end{equation}
Then the solution of~\eqref{MSC_local} with the local model~\eqref{MSC_quadratic_local_model} is
$w_j=R_j\Pi_j^*(f-Av)$.
Consequently, \cref{Alg:PSC,Alg:SSC} reduce to \cref{Alg:MSC_linear_PSC,Alg:MSC_linear_SSC}, respectively.
If $R_j=A_j^{-1}$, then~\eqref{MSC_local} agrees with the exact local problem~\eqref{MSC_local_exact}.
\end{remark}

An identity for PSC for convex optimization that generalizes \cref{Thm:PSC_convergence} was proved in~\cite{LP:2025b,Park:2020}; see \cref{Thm:PSC_convex}.

\begin{theorem}
\label{Thm:PSC_convex}
Suppose that $E\colon V\to\mathbb{R}$ is differentiable and, for every $v\in V$ and $1\leq j\leq J$, the local model $E_j(\cdot;v)\colon V_j\to\mathbb{R}$ is convex and differentiable and satisfies
\begin{equation*}
E_j(0;v)=E(v),
\quad
E_j'(0;v)=\Pi_j^*E'(v).
\end{equation*}
For $w_j\in V_j$, define
\begin{equation*}
d_j(w_j;v)
=E_j(w_j;v)-E(v)-(E'(v),\Pi_jw_j).
\end{equation*}
Then, for $v\in V$, we have
\begin{equation*}
\hat{w}:=\sum_{j=1}^J\Pi_j\hat{w}_j
\in\operatornamewithlimits{\arg\min}_{w\in V}
\left\{
(E'(v),w)
+
\inf_{w=\sum_{j=1}^J\Pi_jw_j}
\sum_{j=1}^Jd_j(w_j;v)
\right\},
\end{equation*}
where $\hat{w}_j$ is the solution of~\eqref{MSC_local} with $u^{\mathrm{old}}=v$, for $1\leq j\leq J$.
Moreover, we have
\begin{equation*}
\inf_{\hat{w}=\sum_{j=1}^J\Pi_jw_j}
\sum_{j=1}^Jd_j(w_j;v)
=
\sum_{j=1}^Jd_j(\hat{w}_j;v).
\end{equation*}
\end{theorem}

Whereas \cref{Thm:PSC_convex} generalizes the linear PSC identity, a sharp convergence estimate for SSC for convex optimization that generalizes the Xu--Zikatanov identity in \cref{Thm:SSC_convergence} remains open.
The convergence of SSC for convex optimization has nevertheless been studied in various settings.
Smooth problems were considered in~\cite{Tai:2003,TE:1998}, and inexact local problems related to the full approximation scheme were studied in~\cite{CHW:2020}.
Constrained problems were treated in~\cite{Badea:2006,BTW:2003}, while composite problems were analyzed in~\cite{Carstensen:1997} and, under more general nonsmooth assumptions, in~\cite{BK:2012}.

The expanded formulation~\eqref{MSC_linear_expanded} extends directly to~\eqref{model_problem_1}.
Define
\begin{equation}
\label{model_problem_1_expanded}
\min_{v_j\in V_j}
\undertilde{E}(\undertilde{v}),
\quad \text{ where } \quad
\undertilde{E}(\undertilde{v})=E(\Pi\undertilde{v}),
\quad \undertilde{v} = (v_1, \dots, v_J) \in \undertilde{V}.
\end{equation}
Problem~\eqref{model_problem_1_expanded} is also an instance of~\eqref{model_problem_1} on the product space $\undertilde{V}$.
Since $\Pi$ is surjective, $u$ solves~\eqref{model_problem_1} if and only if every representation $u=\Pi\undertilde{u}$ solves~\eqref{model_problem_1_expanded}.
Let $\undertilde{E}_j(v_j;\undertilde{v}^{\mathrm{old}})$ be a local block model with a unique minimizer that approximates
$\undertilde{E}(v_1^{\mathrm{old}},\dots,v_{j-1}^{\mathrm{old}},v_j,v_{j+1}^{\mathrm{old}},\dots,v_J^{\mathrm{old}})$.
The expanded Jacobi and expanded Gauss--Seidel methods are presented in \cref{Alg:MSC_block_Jacobi,Alg:MSC_block_GS}, respectively.

\begin{algorithm}
\caption{Expanded Jacobi method for~\eqref{model_problem_1_expanded}}
\label{Alg:MSC_block_Jacobi}
\begin{algorithmic}[]
\STATE Given $\tau>0$ and choose $\undertilde{u}^{(0)}\in\undertilde{V}$.
\FOR{$n=0,1,2,\dots$}
    \FOR{$j=1,2,\dots,J$ \textbf{in parallel}}
        \STATE \vspace{-0.5cm}
        \begin{equation*}
        \hat{u}_j^{(n+1)}
        \in\operatornamewithlimits{\arg\min}_{v_j\in V_j}
        \undertilde{E}_j(v_j;\undertilde{u}^{(n)})
        \end{equation*}
        \vspace{-0.3cm}
    \ENDFOR
    \STATE Set $\hat{\undertilde{u}}^{(n+1)}=(\hat{u}_1^{(n+1)},\dots,\hat{u}_J^{(n+1)})\in\undertilde{V}$.
    \STATE \vspace{-0.5cm}
    \begin{equation*}
    \undertilde{u}^{(n+1)}
    =(1-\tau)\undertilde{u}^{(n)}+\tau\hat{\undertilde{u}}^{(n+1)}
    \end{equation*}
    \vspace{-0.3cm}
\ENDFOR
\end{algorithmic}
\end{algorithm}

\begin{algorithm}
\caption{Expanded Gauss--Seidel method for~\eqref{model_problem_1_expanded}}
\label{Alg:MSC_block_GS}
\begin{algorithmic}[]
\STATE Choose $\undertilde{u}^{(0)}\in\undertilde{V}$.
\FOR{$n=0,1,2,\dots$}
    \FOR{$j=1,2,\dots,J$}
        \STATE \vspace{-0.5cm}
        \begin{equation*}
        u_j^{(n+1)}
        \in\operatornamewithlimits{\arg\min}_{v_j\in V_j}
        \undertilde{E}_j \left( v_j;
        (u_1^{(n+1)},\dots,u_{j-1}^{(n+1)},u_j^{(n)},\dots,u_J^{(n)}) \right)
        \end{equation*}
        \vspace{-0.3cm}
    \ENDFOR
\ENDFOR
\end{algorithmic}
\end{algorithm}

The block Gauss--Seidel method is also called cyclic block coordinate descent and was analyzed in~\cite{BT:2013,Tseng:2001}, while parallel block coordinate descent methods were studied in~\cite{FR:2016,RT:2016}; see also~\cite{Wright:2015} for historical remarks.

The following result extends the linear equivalence between subspace correction and block methods to convex optimization.

\begin{theorem}
\label{Thm:expanded_system}
Suppose that the local models satisfy
\begin{equation}
\label{Thm1:expanded_system}
E_j(w_j;\Pi\undertilde{v}^{\mathrm{old}})
=
\undertilde{E}_j(v_j^{\mathrm{old}}+w_j;\undertilde{v}^{\mathrm{old}}),
\quad w_j\in V_j,\ \undertilde{v}^{\mathrm{old}}\in\undertilde{V}.
\end{equation}
Then PSC~(\cref{Alg:PSC}) is equivalent to the expanded Jacobi method~(\cref{Alg:MSC_block_Jacobi}), and SSC~(\cref{Alg:SSC}) is equivalent to the expanded Gauss--Seidel method~(\cref{Alg:MSC_block_GS}).
More precisely, if $u^{(0)}=\Pi\undertilde{u}^{(0)}$, then the corresponding iterates satisfy
\begin{equation*}
u^{(n)}=\Pi\undertilde{u}^{(n)},
\quad n\geq1.
\end{equation*}
\end{theorem}
\begin{proof}
For the expanded Jacobi method, suppose that $u^{(n)}=\Pi\undertilde{u}^{(n)}$.
By~\eqref{Thm1:expanded_system} and uniqueness of the local minimizers, we have
$w_j^{(n+1)}=\hat{u}_j^{(n+1)}-u_j^{(n)}$.
Consequently,
\begin{equation*}
u^{(n+1)}
=u^{(n)}+\tau\sum_{j=1}^J\Pi_jw_j^{(n+1)}
=(1-\tau)\Pi\undertilde{u}^{(n)}
+\tau\Pi\hat{\undertilde{u}}^{(n+1)}
=\Pi\undertilde{u}^{(n+1)}.
\end{equation*}
For the expanded Gauss--Seidel method, the same argument at the $j$th local update gives
$w_j^{(n+1)}=u_j^{(n+1)}-u_j^{(n)}$.
Induction over $j$ therefore yields
$u^{(n+\frac{j}{J})}=\sum_{i=1}^j\Pi_i u_i^{(n+1)}+\sum_{i=j+1}^J\Pi_i u_i^{(n)}$, which proves the result.
\end{proof}

The expanded problem~\eqref{model_problem_1_expanded} may be semicoercive even if~\eqref{model_problem_1} is coercive, due to the possible nontrivial kernel of $\Pi$.
This connects the convex theory to singular linear systems~\cite{LWXZ:2008,WLXZ:2008}; recent convergence results for semicoercive convex optimization can be found in~\cite{LP:2025b}.

\begin{remark}
\label{Rem:acceleration}
The expanded system viewpoint also permits techniques developed for block methods to be transferred to subspace correction methods.
Randomized subspace correction methods were considered in~\cite{GO:2012,HXZ:2019}.
Randomized block coordinate methods were introduced in~\cite{Nesterov:2012,RT:2014}, and parallel variants were studied in~\cite{FR:2015,FR:2016,RT:2016}.
A general framework for randomized iterative methods for linear systems, including randomized Kaczmarz and coordinate descent methods, was developed in~\cite{GR:2015}.
Acceleration of cyclic block methods was considered in~\cite{BT:2013,CP:2015,ST:2016}, while acceleration of parallel methods was studied in~\cite{LP:2019,PX:2023}.
A broader upper-bound minimization framework was proposed in~\cite{HWRL:2017,RHL:2013}.
Recent work on acceleration of general subspace correction methods can be found in~\cite{Park:2021,Park:2022a}.
\end{remark}

\section{Operator splitting methods}
\label{Sec:Splitting}
Operator splitting methods, originating from the pioneering Peaceman--Rachford and Douglas--Rachford methods~\cite{DR:1956,PR:1955}, have been widely applied not only to linear problems but also to more general settings in convex optimization and monotone operator theory~\cite{CKCH:2023,EB:1992,LM:1979}.
Alternating projection methods arise as special cases of operator splitting methods when the convex functions are indicator functions of convex sets~\cite{CP:2009}.
In this section, we show that, under suitable dual formulations, these methods arise as dualizations of certain subspace correction or block methods.
This interpretation also suggests multi-block variants and a parallel variant obtained by dualizing corresponding successive and parallel subspace correction methods.

\subsection{Optimization of the sum of multiple convex functions}
\label{Subsec:multiple_convex}
We begin with an optimization problem involving the sum of multiple convex functions.
The extension of splitting algorithms from two convex functions to multiple convex functions is regarded as a nontrivial task~\cite{LM:1979,Ryu:2020}.

Consider the convex optimization problem
\begin{equation}
\label{multiple_convex_primal}
\min_{v \in V} \left\{ F(v) + \sum_{j=1}^J G_j (B_j v) \right\}.
\end{equation}
Here, each $W_j$ is a real Hilbert space, $B_j\colon V\to W_j$ is a bounded linear operator, and $F\colon V\to\overline{\mathbb{R}}$ and $G_j\colon W_j\to\overline{\mathbb{R}}$ are proper, convex, and lower semicontinuous.
Problem~\eqref{multiple_convex_primal} is a multi-term specialization of~\eqref{model_problem_2}, obtained by setting
\begin{equation*}
W=\undertilde{W}:=\prod_{j=1}^J W_j,
\quad
G(\undertilde{w})=\sum_{j=1}^J G_j(w_j),
\quad
Bv=(B_1v,\dots,B_Jv),
\end{equation*}
where $\undertilde{w}=(w_1,\dots,w_J)\in\undertilde{W}$.
Under the corresponding qualification condition, its dual problem is given by
\begin{equation}
\label{multiple_convex_dual}
\min_{q_j \in W_j^*} \left\{ F^* \left( - \sum_{j=1}^J B_j^* q_j \right) + \sum_{j=1}^J G_j^* (q_j) \right\},
\end{equation}
where $\undertilde{W}^* = \prod_{j=1}^J W_j^*$ and
$\undertilde{q}=(q_1,\dots,q_J)\in\undertilde{W}^*$.

For the algorithms below, we assume in addition that \(F\) is real-valued, strongly convex, and smooth, ensuring that~\eqref{Legendre} holds, i.e., $F'$ is invertible with inverse $(F^*)'$.
We also assume that the subproblems involving $F$ and a single $G_j\circ B_j$ admit closed-form solutions or can be solved efficiently, making~\eqref{multiple_convex_primal} suitable for operator splitting methods.
A recent review of splitting algorithms for solving~\eqref{multiple_convex_primal} can be found in~\cite{CKCH:2023}.

An important special case of~\eqref{multiple_convex_primal}, and hence of~\eqref{model_problem_2}, is given by  
\begin{equation}
\label{multiple_convex_special}
\min_{v \in V} \left\{ \sum_{j=1}^J G_j (v) + \frac{\alpha}{2} \| v \|^2 \right\},
\end{equation}
where $\alpha$ is a positive regularization parameter.
Problem~\eqref{multiple_convex_special} frequently appears in machine learning and statistics~\cite{SB:2014}.
A variety of algorithms have been developed to solve~\eqref{multiple_convex_special}, including stochastic gradient methods~\cite{Bertsekas:2011,Bertsekas:2015}.

The corresponding primal--dual relations~\eqref{primal_dual_relations} are written as
\begin{equation}
\label{multiple_primal_dual_relations}
- \sum_{j=1}^J B_j^* p_{j} = F'(u), \quad
B_j u \in \partial G_j^* (p_{j}),
\quad 1 \leq j \leq J.
\end{equation}

We begin with SSC for solving the dual problem~\eqref{multiple_convex_dual}.
We apply SSC to the coordinate blocks of the dual variable space $\undertilde{W}^*=\prod_{j=1}^J W_j^*$.
The $j$th local problem is written as
\begin{equation}
\label{multiple_convex_GS}
    p_j^{(n+1)} \in \operatornamewithlimits{\arg\min}_{q_j \in W_j^*} \left\{ F^* \left( -B_j^* q_j - \sum_{i = 1}^{j-1} B_i^* p_i^{(n+1)} - \sum_{i = j+1}^J B_i^* p_i^{(n)} \right) + G_j^* (q_j) \right\},
\end{equation}
where $\undertilde{p}^{(0)}=(p_1^{(0)},\dots,p_J^{(0)})\in \undertilde{W}^*$ is an initial guess.

Dualizing the local problem~\eqref{multiple_convex_GS} gives the primal subproblem
\begin{equation}
\label{multiple_convex_PR_2}
\min_{v\in V}
\left\{
D_F(v;u^{(n+\frac{j-1}{J})})
-
( B_j^* p_j^{(n)},v)
+
G_j(B_j v)
\right\}.
\end{equation}
Indeed, if
\[
u^{(n+\frac{j}{J})}
=
(F^*)'
\left(
-\sum_{i=1}^j B_i^* p_i^{(n+1)}
-\sum_{i=j+1}^J B_i^* p_i^{(n)}
\right),
\]
then the primal--dual relations for the local pair imply that the solution of~\eqref{multiple_convex_PR_2} is precisely \(u^{(n+\frac{j}{J})}\).
Thus, after setting \(r_j^{(n)}=-B_j^* p_j^{(n)}\), we obtain \cref{Alg:multiple_convex_PR}, which we call the dualization-based Peaceman--Rachford splitting algorithm.

\begin{algorithm}
\caption{Dualization-based Peaceman--Rachford splitting algorithm for~\eqref{multiple_convex_primal}}
\label{Alg:multiple_convex_PR}
\begin{algorithmic}[]
\begin{subequations}
\STATE Choose $u^{(0)} \in V$ and $\undertilde{r}^{(0)}=(r_1^{(0)},\dots,r_J^{(0)}) \in (V^*)^J$.
\FOR{$n=0,1,2,\dots$}
    \FOR{$j=1,2,\dots,J$}
        \STATE \vspace{-0.5cm} \begin{align}
        \label{Alg1:multiple_convex_PR}
        &u^{(n+\frac{j}{J})} = \operatornamewithlimits{\arg\min}_{v \in V} \left\{ D_F (v; u^{(n+\frac{j-1}{J})}) + ( r_j^{(n)}, v ) + G_j (B_j v) \right\} \\
        &r_j^{(n+1)} = r_j^{(n)} + F'( u^{(n+\frac{j}{J})} ) - F'( u^{(n+\frac{j-1}{J})} )
        \end{align}
        \vspace{-0.3cm}
    \ENDFOR
\ENDFOR
\end{subequations}
\end{algorithmic}
\end{algorithm}

By construction, we have the following theorem, which summarizes the dualization relation between SSC and \cref{Alg:multiple_convex_PR}.

\begin{theorem}
\label{Thm:multiple_convex_PR}
The dualization-based Peaceman--Rachford splitting algorithm~(\cref{Alg:multiple_convex_PR}) is a dualization of SSC~\eqref{multiple_convex_GS} for solving~\eqref{multiple_convex_dual}.
More precisely, if
\begin{equation*}
u^{(0)} = (F^*)' \left( - \sum_{j=1}^J B_j^* p_j^{(0)} \right), \quad
r_j^{(0)} = - B_j^* p_j^{(0)} , \quad 1 \leq j \leq J,
\end{equation*}
then we have
\begin{equation*}
u^{(n)} = (F^*)' \left( - \sum_{j=1}^J B_j^* p_j^{(n)} \right), \quad
r_j^{(n)} = - B_j^* p_j^{(n)} ,
\quad 1 \leq j \leq J,\  n \geq 1.
\end{equation*}
\end{theorem}

\begin{remark}
\label{Rem:conventional_PR}
When \cref{Alg:multiple_convex_PR} is applied to~\eqref{multiple_convex_special} with \(J=2\), \(\alpha=2\), and \(B_1=B_2=I\), it reduces to the conventional Peaceman--Rachford splitting algorithm~\cite{CKCH:2023,LM:1979}.
Indeed, writing
\[
\bar{G}_j(v)=G_j(v)+\frac12\|v\|^2,\quad j=1,2,
\]
problem~\eqref{multiple_convex_special} becomes
\begin{equation}
\label{two_convex}
\min_{v\in V}\{\bar{G}_1(v)+\bar{G}_2(v)\}.
\end{equation}
Problem~\eqref{two_convex} is an instance of~\eqref{model_problem_2}.
Under the initialization
\begin{equation}
\label{PR_conventional_initialization}
u^{(0)}=\frac12(r_1^{(0)}+r_2^{(0)}),
\end{equation}
one has \(u^{(n)}=\frac12(r_1^{(n)}+r_2^{(n)})\). Setting
$z^{(n)}=2u^{(n)}-r_1^{(n)}$,
the iteration in \cref{Alg:multiple_convex_PR} takes the following form:
\begin{subequations}
\label{PR_conventional}
\begin{align}
 u^{(n+\frac12)}
 &=
 \operatornamewithlimits{\arg\min}_{v\in V}
 \left\{\frac12\|v-z^{(n)}\|^2+\bar G_1(v)\right\}, \\
 u^{(n+1)}
 &=
 \operatornamewithlimits{\arg\min}_{v\in V}
 \left\{\frac12\|v-(2u^{(n+\frac12)}-z^{(n)})\|^2+\bar G_2(v)\right\}, \\
 z^{(n+1)}
 &=
 2u^{(n+1)}-2u^{(n+\frac12)}+z^{(n)}.
\end{align}
\end{subequations}
This is equivalent to the conventional Peaceman--Rachford splitting formulation~(see, e.g.,~\cite[Algorithm~I]{LM:1979}).
\end{remark}

\subsection{Relaxed splitting algorithms: Douglas--Rachford splitting}
\label{Subsec:relaxed-DR}
In the conventional Peaceman--Rachford splitting algorithm~\eqref{PR_conventional} for solving~\eqref{two_convex}, introducing a relaxation step yields the Douglas--Rachford splitting algorithm (see, e.g.,~\cite[Algorithm~II]{LM:1979} and~\cite[Equation~(5.22)]{CKCH:2023}):  
\begin{subequations}
\label{DR_conventional}
\begin{align}
 \tilde{u}^{(n+\frac12)}
 &=
 \operatornamewithlimits{\arg\min}_{v\in V}
 \left\{
 \frac12\|v-z^{(n)}\|^2+\bar G_1(v)
 \right\}, \\
 \tilde{u}^{(n+1)}
 &=
 \operatornamewithlimits{\arg\min}_{v\in V}
 \left\{
 \frac12\|v-(2\tilde{u}^{(n+\frac12)}-z^{(n)})\|^2
 +\bar G_2(v)
 \right\}, \\
 \hat z^{(n+1)}
 &=
 2\tilde{u}^{(n+1)}
 -2\tilde{u}^{(n+\frac12)}
 +z^{(n)}, \\
 z^{(n+1)}
 &=
 (1-\tau)z^{(n)}
 +\tau\hat z^{(n+1)},
\end{align}
\end{subequations}
where $0<\tau<1$ is a relaxation parameter.  
This algorithm is widely used due to its robust convergence properties under relatively weak conditions~\cite{Combettes:2009}.  
Moreover, the connections among the Douglas--Rachford splitting algorithm, the alternating direction method of multipliers~\cite{Glowinski:2014,WYZ:2019}, and the Chambolle--Pock primal--dual algorithm~\cite{CP:2011,CP:2016b} have been well studied; see~\cite{CP:2011,CKCH:2023}.  

We next show that a Douglas--Rachford variant of \cref{Alg:multiple_convex_PR}, called the \emph{dualization-based Douglas--Rachford splitting algorithm}, has an analogous interpretation as a dualization of SSC for~\eqref{multiple_convex_dual}.

\begin{algorithm}
\caption{Dualization-based Douglas--Rachford splitting algorithm for~\eqref{multiple_convex_primal}}
\label{Alg:multiple_convex_DR}
\begin{algorithmic}[]
\begin{subequations}
\STATE Given $0<\tau<1$:
\STATE Choose $u^{(0)} \in V$ and $\undertilde{r}^{(0)}=(r_1^{(0)},\dots,r_J^{(0)}) \in (V^*)^J$.
\FOR{$n=0,1,2,\dots$}
    \STATE \vspace{-0.5cm} \begin{equation}
    \label{Alg1:multiple_convex_DR}
    \hat{u}^{(n+1, 0)} = u^{(n)}
    \end{equation}
    \vspace{-0.3cm}
    \FOR{$j=1,2,\dots,J$}
        \STATE \vspace{-0.5cm} \begin{align}
        &\hat{u}^{(n+1, j)} = \operatornamewithlimits{\arg\min}_{v \in V} \left\{ D_F (v; \hat{u}^{(n+1, j-1)}) + ( r_j^{(n)}, v ) + G_j (B_j v) \right\}\\
        \label{Alg3:multiple_convex_DR}
        &r_j^{(n+1)} = r_j^{(n)} + \tau ( F'( \hat{u}^{(n+1, j)} ) - F'( \hat{u}^{(n+1, j-1)} )  )
        \end{align}
        \vspace{-0.3cm}
    \ENDFOR
    \STATE \vspace{-0.5cm} \begin{align}
    \label{Alg4:multiple_convex_DR}
    &u^{(n+1)} = (F^*)' \left( \sum_{j=1}^J r_j^{(n+1)} \right)
    \end{align}
    \vspace{-0.3cm}
\ENDFOR
\end{subequations}
\end{algorithmic}
\end{algorithm}

To separate the unrelaxed splitting step from the relaxation step in
\Cref{Alg:multiple_convex_DR}, we define
\begin{equation}
\label{multiple_convex_DR_r_hat1}
\hat{r}_j^{(n+1)}
=
r_j^{(n)}
+
F'( \hat{u}^{(n+1, j)} )
-
F'( \hat{u}^{(n+1, j-1)} ),
\quad
1 \leq j \leq J,\  n \geq 0.
\end{equation}
Then the update~\eqref{Alg3:multiple_convex_DR} can be written as
\begin{equation}
\label{multiple_convex_DR_r_hat2}
r_j^{(n+1)}
=
(1-\tau)r_j^{(n)}
+
\tau \hat{r}_j^{(n+1)},
\quad
1 \leq j \leq J,\  n \geq 0.
\end{equation}

\begin{remark}
\label{Rem:conventional_DR}
When \cref{Alg:multiple_convex_DR} is applied to~\eqref{multiple_convex_special} with \(J=2\), \(\alpha=2\), and \(B_1=B_2=I\), the same change of variables as in \cref{Rem:conventional_PR} shows that it reduces to the conventional Douglas--Rachford splitting algorithm~\eqref{DR_conventional}.
Indeed, under the initialization~\eqref{PR_conventional_initialization}, set
\[
z^{(n)}=r_2^{(n)},
\quad
\hat z^{(n+1)}=\hat r_2^{(n+1)}.
\]
Then
\[
\tilde{u}^{(n+\frac12)}
=
\hat u^{(n+1,1)},
\quad
\tilde{u}^{(n+1)}
=
\hat u^{(n+1,2)}
\]
satisfy the two minimization substeps in~\eqref{DR_conventional}, and
\[
\hat z^{(n+1)}
=
2\tilde{u}^{(n+1)}
-
2\tilde{u}^{(n+\frac12)}
+
z^{(n)}.
\]
Moreover,~\eqref{multiple_convex_DR_r_hat2} gives
\[
z^{(n+1)}
=
(1-\tau)z^{(n)}
+
\tau\hat z^{(n+1)}.
\]
Thus, \(z^{(n)}\) satisfies the iteration~\eqref{DR_conventional}.
\end{remark}

We now introduce the dual algorithm whose dualization gives
\Cref{Alg:multiple_convex_DR}.
Let \(\hat{\undertilde{p}}^{(n+1)}=(\hat p_1^{(n+1)},\dots,\hat p_J^{(n+1)})\in\undertilde{W}^*\) be the intermediate iterate obtained by one SSC sweep~\eqref{multiple_convex_GS} applied to the dual problem~\eqref{multiple_convex_dual}.
The relaxed SSC is then defined by
\begin{equation}
\label{multiple_convex_RGS}
\undertilde{p}^{(n+1)}
=
(1-\tau)\undertilde{p}^{(n)}
+
\tau\hat{\undertilde{p}}^{(n+1)}.
\end{equation}

We next establish that \cref{Alg:multiple_convex_DR} is a dualization of the relaxed SSC defined by
\eqref{multiple_convex_GS} and~\eqref{multiple_convex_RGS}.

\begin{theorem}
\label{Thm:multiple_convex_DR}
The dualization-based Douglas--Rachford splitting algorithm~(\cref{Alg:multiple_convex_DR}) is a dualization of the relaxed SSC~(see~\eqref{multiple_convex_GS} and~\eqref{multiple_convex_RGS}) for solving~\eqref{multiple_convex_dual}.
More precisely, if
\begin{equation}
\label{multiple_convex_DR_initial}
u^{(0)} = (F^*)' \left( - \sum_{j=1}^J B_j^*p_j^{(0)} \right), \quad
r_j^{(0)} = - B_j^* p_j^{(0)} ,
\quad 1 \leq j \leq J,
\end{equation}
then we have
\begin{equation}
\label{Thm1:multiple_convex_DR}
u^{(n)} = (F^*)' \left( - \sum_{j=1}^J B_j^*p_j^{(n)} \right), \quad
r_j^{(n)} = - B_j^* p_j^{(n)} ,
\quad 1 \leq j \leq J,\  n \geq 1.
\end{equation}
\end{theorem}
\begin{proof}
Since the proof is similar to that of \cref{Thm:multiple_convex_PR}, we provide only a sketch.  
Assuming that~\eqref{Thm1:multiple_convex_DR} holds for some $n \geq 0$, we can establish the following by applying an argument analogous to the proof of \cref{Thm:multiple_convex_PR}:  
\begin{equation*}
\hat{u}^{(n+1,j)} = (F^*)' \left( - \sum_{i=1}^j B_i^* \hat{p}_i^{(n+1)} - \sum_{i=j+1}^J B_i^* p_i^{(n)} \right),\ 
\hat{r}_j^{(n+1)} = -B_j^* \hat{p}_j^{(n+1)},\ 1 \leq j \leq J,
\end{equation*}
where $\hat{r}_j^{(n+1)}$ was defined in~\eqref{multiple_convex_DR_r_hat1}.
Then, the validity of~\eqref{Thm1:multiple_convex_DR} for $n+1$ follows directly from~\eqref{Alg4:multiple_convex_DR},~\eqref{multiple_convex_DR_r_hat2}, and~\eqref{multiple_convex_RGS}.  
\end{proof}

\begin{remark}
\label{Rem:multiple_convex_DR}
When applied to~\eqref{multiple_convex_special}, i.e., when $ F(v) = \frac{\alpha}{2} \| v \|^2 $, the step~\eqref{Alg4:multiple_convex_DR} can be replaced with a simpler update rule:
\begin{equation*}
u^{(n+1)} = (1 - \tau) u^{(n)} + \tau \hat{u}^{(n+1,J)},
\end{equation*}
under the assumption~\eqref{multiple_convex_DR_initial}.  
To justify this simplification, suppose that~\eqref{Thm1:multiple_convex_DR} holds for $ n $. Then, we derive
\begin{multline*}
u^{(n+1)} \stackrel{\eqref{Alg4:multiple_convex_DR}}{=} \frac{1}{\alpha} \sum_{j=1}^J r_j^{(n+1)}
\stackrel{\eqref{Alg3:multiple_convex_DR}}{=}
\frac{1}{\alpha} \sum_{j=1}^J r_j^{(n)} + \tau (\hat{u}^{(n+1,J)} - \hat{u}^{(n+1,0)}) \\
\stackrel{\eqref{Thm1:multiple_convex_DR}}{=} u^{(n)} + \tau (\hat{u}^{(n+1,J)} - \hat{u}^{(n+1,0)})
\stackrel{\eqref{Alg1:multiple_convex_DR}}{=} (1 - \tau) u^{(n)} + \tau \hat{u}^{(n+1,J)}.
\end{multline*}
\end{remark}

\begin{remark}
\label{Rem:DR_extension}
There is a substantial literature on generalizations of Douglas--Rachford splitting.
Examples include projection variants for feasibility problems~\cite{DP:2018}, three-operator extensions~\cite{DY:2017}, and extensions for composite and parallel-sum-type monotone operators~\cite{BH:2012}.
\end{remark}

\subsection{New parallel operator splitting algorithm}
The dualization viewpoint also suggests a parallel variant of the splitting algorithms above.
Instead of applying SSC to the dual problem~\eqref{multiple_convex_dual}, we apply PSC~(\cref{Alg:PSC}) with the coordinate-block decomposition of $\undertilde{W}^*=\prod_{j=1}^J W_j^*$.
The resulting local problems and relaxation step are given as follows:
\begin{subequations}
\label{multiple_convex_Jacobi}
\begin{align}
    &\hat{p}_j^{(n+1)} \in
    \operatornamewithlimits{\arg\min}_{q_j \in W_j^*}
    \left\{
        F^*
        \left(
            - B_j^* q_j
            -
            \sum_{i \neq j} B_i^* p_i^{(n)}
        \right)
        +
        G_j^* (q_j)
    \right\},
    \quad 1\leq j\leq J, \\
    &\displaystyle
    \undertilde{p}^{(n+1)}
    =
    (1-\tau)\undertilde{p}^{(n)}
    +
    \tau \hat{\undertilde{p}}^{(n+1)},
\end{align}
\end{subequations}
where $\undertilde{p}^{(0)}=(p_1^{(0)},\dots,p_J^{(0)})\in\undertilde{W}^*$ is an initial guess and
$\hat{\undertilde{p}}^{(n+1)}=(\hat p_1^{(n+1)},\dots,\hat p_J^{(n+1)})\in\undertilde{W}^*$.
Dualizing~\eqref{multiple_convex_Jacobi} yields the following parallel splitting algorithm.
We then establish the corresponding dualization relation.

\begin{algorithm}
\caption{Parallel Douglas--Rachford splitting algorithm for~\eqref{multiple_convex_primal}}
\label{Alg:multiple_convex_parallel}
\begin{algorithmic}[]
\begin{subequations}
\STATE Given $\tau > 0$:
\STATE Choose $u^{(0)} \in V$ and $\undertilde{r}^{(0)}=(r_1^{(0)},\dots,r_J^{(0)}) \in (V^*)^J$.
\FOR{$n=0,1,2,\dots$}
    \FOR{$j=1,2,\dots,J$ \textbf{in parallel}}
        \STATE \vspace{-0.5cm} \begin{align}
            \label{Alg1:multiple_convex_parallel}
            &\hat{u}_j^{(n+1)} = \operatornamewithlimits{\arg\min}_{v \in V} \left\{ D_F (v; u^{(n)}) + ( r_j^{(n)}, v ) + G_j (B_j v) \right\} \\
            &r_j^{(n+1)} = r_j^{(n)} + \tau ( F'( \hat{u}_j^{(n+1)}) - F' (u^{(n)}) )
        \end{align}
        \vspace{-0.3cm}
    \ENDFOR
    \STATE \vspace{-0.5cm} \begin{equation*}
        u^{(n+1)} = (F^*)' \left( \sum_{j=1}^J r_j^{(n+1)} \right)
    \end{equation*}
    \vspace{-0.3cm}
\ENDFOR
\end{subequations}
\end{algorithmic}
\end{algorithm}

\begin{theorem}
\label{Thm:multiple_convex_parallel}
The parallel Douglas--Rachford splitting algorithm~(\cref{Alg:multiple_convex_parallel}) is a dualization of PSC~\eqref{multiple_convex_Jacobi} for solving~\eqref{multiple_convex_dual}.
More precisely, if
\begin{equation}
\label{multiple_convex_parallel_initial}
u^{(0)} = (F^*)' \left( - \sum_{j=1}^J B_j^*p_j^{(0)} \right), \quad
r_j^{(0)} = - B_j^* p_j^{(0)} ,
\quad 1 \leq j \leq J,
\end{equation}
then we have
\begin{equation*}
u^{(n)} = (F^*)' \left( - \sum_{j=1}^J B_j^*p_j^{(n)} \right), \quad
r_j^{(n)} = - B_j^* p_j^{(n)} ,
\quad 1 \leq j \leq J,\ n \geq 1.
\end{equation*}
\end{theorem}
\begin{proof}
The proof is analogous to that of \cref{Thm:multiple_convex_PR}, with the local problems solved in parallel, so we omit the details.
\end{proof}

\begin{remark}
\label{Rem:multiple_convex_parallel}
As in \cref{Rem:multiple_convex_DR}, when \cref{Alg:multiple_convex_parallel} is applied to~\eqref{multiple_convex_special}, the update for $ u^{(n+1)} $ simplifies to
\begin{equation*}
u^{(n+1)} = (1 - \tau J) u^{(n)} + \tau \sum_{j=1}^J \hat{u}_j^{(n+1)},
\end{equation*}
under the assumption~\eqref{multiple_convex_parallel_initial}.
\end{remark}

In Section~9 of a recent survey~\cite{CKCH:2023} on splitting algorithms for convex optimization, two techniques for designing parallel algorithms were introduced. Both techniques are based on reformulations of the problem~\eqref{multiple_convex_primal}, which involves $J$ convex functions, into simpler structures. Specifically, two formulations are considered.
The first is given by
\begin{equation}
\label{parallel_reformulation_1}
\min_{v_j\in V}
\left\{
F(v_1)
+
G(\undertilde{B}\undertilde{v})
+
\Psi(\undertilde{v})
\right\},
\end{equation}
where $\undertilde{W} = \prod_{j=1}^J W_j$.
The functions
$G \colon \undertilde{W} \to \overline{\mathbb{R}}$ and
$\Psi \colon V^J \to \overline{\mathbb{R}}$,
together with the linear operator
$\undertilde{B} \colon V^J \to \undertilde{W}$, are defined by
\begin{multline*}
G(\undertilde{w}) = \sum_{j=1}^J G_j(w_j), \quad
\Psi(\undertilde{v}) =
\begin{cases}
0, & \textrm{if } v_1 = \dots = v_J, \\
\infty, & \textrm{otherwise},
\end{cases}
\quad
\undertilde{B} = \operatorname{diag}(B_j), \\
\undertilde{v}=(v_1,\dots,v_J)\in V^J,\quad
\undertilde{w}=(w_1,\dots,w_J)\in\undertilde{W}.
\end{multline*}
Problem~\eqref{parallel_reformulation_1} is an instance of~\eqref{model_problem_2} on the product space $V^J$.
The second formulation is given by
\begin{equation}
\label{parallel_reformulation_2}
\min_{v\in V}
\left\{
F(v)+(G\circ\underline{B})(v)
\right\},
\quad
\underline{B}\colon V\to\undertilde{W},
\quad
\underline{B}v=(B_1v,\dots,B_Jv).
\end{equation}
Problem~\eqref{parallel_reformulation_2} is an instance of~\eqref{model_problem_2}.
Applying well-known splitting algorithms for two functions to~\eqref{parallel_reformulation_1} or~\eqref{parallel_reformulation_2} leads to parallel splitting algorithms; see~\cite[Section~9]{CKCH:2023} for details.
While designing sequential splitting algorithms from these reformulations is not straightforward, the dualization approach introduced in this paper provides a common perspective on constructing both parallel and sequential algorithms.

\begin{remark}
\label{Rem:multiple_primal_dual_splitting}
Problem~\eqref{multiple_convex_primal} is also a standard model in primal--dual splitting frameworks, such as those of Esser~\cite{EZC:2010}, Chambolle--Pock~\cite{CP:2011}, Condat--V{\~u}~\cite{Condat:2013,Vu:2011}, and Combettes--Pesquet~\cite{CPesquet:2011}.
In general, the iterates of these primal--dual methods need not satisfy the first relation in~\eqref{multiple_primal_dual_relations} before convergence.
By contrast, in the dualization-based methods introduced in this section, that relation is imposed at every iterate.
\end{remark}

\subsection{Special cases}
The general splitting algorithms above include several classical methods as special cases.
We first consider linear systems with multiple operators and then alternating projection methods.

\subsubsection{Linear systems}
\label{Subsec:linear_splitting}
We consider the linear system
\begin{equation}
\label{multiple_linear_primal}
\left(\sum_{j=1}^J A_j+\alpha I\right)u=f,
\end{equation}
where each $A_j\colon V\to V^*$ is SPD, $\alpha>0$, and $f\in V^*$.
The linear system~\eqref{multiple_linear_primal} is the optimality condition for a quadratic instance of~\eqref{model_problem_2}.
In~\cite{DKV:1963}, the following multi-block Peaceman--Rachford splitting algorithm was considered for solving~\eqref{multiple_linear_primal}.
Choose $u^{(-1+\frac{j}{J})}\in V$ for $1\leq j\leq J$ and set

\begin{equation}
\label{Alg:multiple_linear_PR}
u^{(n+\frac{j}{J})}
=
(A_j+\alpha I)^{-1}
\left(
f-
\sum_{i=1}^{j-1}A_i u^{(n+\frac{i}{J})}
-
\sum_{i=j+1}^{J}A_i u^{(n-1+\frac{i}{J})}
\right),
\quad n\geq0,\ 1\leq j\leq J.
\end{equation}

As a special case of \cref{Thm:multiple_convex_PR}, we deduce that iteration~\eqref{Alg:multiple_linear_PR} is a dualization of SSC for the dual system
\begin{equation}
\label{multiple_linear_dual}
\alpha A_j^{-1} p_j + \sum_{i=1}^J p_i = f,
\quad 1 \leq j \leq J.
\end{equation}
System~\eqref{multiple_linear_dual} is the optimality system for a quadratic instance of~\eqref{model_problem_1} on the product space $(V^*)^J$.
Its SSC update is given by
\begin{equation*}
p_j^{(n+1)}
=
(\alpha A_j^{-1}+I)^{-1}
\left(
f-
\sum_{i=1}^{j-1}p_i^{(n+1)}
-
\sum_{i=j+1}^Jp_i^{(n)}
\right),
\quad n \geq 0,\ 1 \leq j \leq J.
\end{equation*}
With $p_j^{(0)}=A_j u^{(-1+\frac{j}{J})}$, the iterates satisfy
\begin{equation*}
u^{(n)}
=
\frac{1}{\alpha}
\left(
f-\sum_{j=1}^Jp_j^{(n)}
\right),
\quad n\geq1,
\end{equation*}
which is the defining dualization relation for this linear case.

\begin{remark}
\label{Rem:multiple_convex_PR}
The dualization-based Peaceman--Rachford splitting algorithm~(\cref{Alg:multiple_convex_PR}) generalizes the linear iteration~\eqref{Alg:multiple_linear_PR}.
Indeed, for~\eqref{multiple_linear_primal}, set
\begin{equation*}
F(v)=\frac{\alpha}{2}\|v\|^2-(f,v),
\quad
G_j(v)=\frac12(A_jv,v),
\quad
B_j=I.
\end{equation*}
Then $G_j'(v)=A_jv$, and the second primal--dual relation gives
\begin{equation*}
r_j^{(n)}=-A_j u^{(n-1+\frac{j}{J})},
\quad n\geq1.
\end{equation*}
Substituting this identity into~\eqref{Alg1:multiple_convex_PR} and using the iterate relation in \cref{Thm:multiple_convex_PR} yields~\eqref{Alg:multiple_linear_PR}.
More generally, if each $G_j$ is differentiable, the auxiliary variable $r_j^{(n)}$ can be eliminated by the same relation.
\end{remark}

\begin{remark}
\label{Rem:multiple_linear_PR}
In~\cite{DKV:1963}, iteration~\eqref{Alg:multiple_linear_PR} was derived as a block Gauss--Seidel method for a nonsymmetric block system.
When $J=2$, it reduces to the classical Peaceman--Rachford algorithm for two operators~\cite{BVY:1962,PR:1955}.
Thus, the relation above also connects SSC for that formulation with SSC for the symmetric dual system~\eqref{multiple_linear_dual}.
\end{remark}

\subsubsection{Alternating projection methods}
\label{Subsec:projection_special}
Projection onto convex sets is a fundamental problem with widespread applications.
The von Neumann algorithm~\cite{Neumann:1950}, introduced in~\cref{Alg:MSC_von_Neumann} for subspaces, is the classical alternating projection method.
The Dykstra algorithm~\cite{BD:1986,Dykstra:1983} extends cyclic projection to general closed convex sets by using correction terms, and its relation with block coordinate descent was studied in~\cite{GM:1989,Han:1988,Tibshirani:2017}.
A general account of alternating projection methods is given in~\cite{BB:1996}.
Related parallel algorithms were considered in~\cite{CP:2009,Tibshirani:2017}.

Let $K_j$, $1\leq j\leq J$, be closed convex subsets of $V$ with nonempty intersection.
The corresponding projection problem is given by
\begin{equation}
\label{POCS}
\min_{v\in V}
\left\{
\frac12\|v-f\|^2
+
\chi_{\bigcap_{j=1}^JK_j}(v)
\right\}.
\end{equation}
Problem~\eqref{POCS} is a special case of~\eqref{multiple_convex_primal}, and hence of~\eqref{model_problem_2}, with
$W_j=V$, $F(v)=\frac12\|v-f\|^2$, $G_j=\chi_{K_j}$, and $B_j=I$~\cite{CP:2009}.
After omitting an additive constant, the dual problem~\eqref{multiple_convex_dual} becomes
\begin{equation}
\label{POCS_subsets_dual}
\min_{q_j\in V^*}
\left\{
\frac12\left\|\sum_{j=1}^Jq_j-f\right\|^2
+
\sum_{j=1}^J\chi_{K_j}^*(q_j)
\right\}.
\end{equation}
Problem~\eqref{POCS_subsets_dual} is also an instance of~\eqref{model_problem_2}.
Since $D_F(v;w)=\frac12\|v-w\|^2$ and $F'(v)=v-f$, setting $s_j^{(n)}=-r_j^{(n)}$ in \cref{Alg:multiple_convex_PR} gives the Dykstra algorithm in \cref{Alg:Dykstra}.

\begin{algorithm}
\caption{Dykstra algorithm for~\eqref{POCS}}
\label{Alg:Dykstra}
\begin{algorithmic}[]
\STATE Set $u^{(0)}=f$ and $s_j^{(0)}=0$, $1\leq j\leq J$.
\FOR{$n=0,1,2,\dots$}
    \FOR{$j=1,2,\dots,J$}
        \STATE \vspace{-0.5cm}
        \begin{align*}
        u^{(n+\frac{j}{J})}
        &=
        \operatorname{proj}_{K_j}
        (u^{(n+\frac{j-1}{J})}+s_j^{(n)}),\\
        s_j^{(n+1)}
        &=
        s_j^{(n)}+u^{(n+\frac{j-1}{J})}-u^{(n+\frac{j}{J})}.
        \end{align*}
        \vspace{-0.3cm}
    \ENDFOR
\ENDFOR
\end{algorithmic}
\end{algorithm}
The specialization of \cref{Thm:multiple_convex_PR} shows that \cref{Alg:Dykstra} is a dualization of SSC applied to~\eqref{POCS_subsets_dual} with $\undertilde{p}^{(0)}=0$.
More precisely, we have
\begin{equation*}
u^{(n+\frac{j}{J})}
=
f-
\left(
\sum_{i\leq j}p_i^{(n+1)}
+
\sum_{i>j}p_i^{(n)}
\right),
\quad
s_i^{(n)}=p_i^{(n)},
\quad 1\leq i\leq J.
\end{equation*}
This recovers the equivalence between the Dykstra algorithm and block coordinate descent established in~\cite{GM:1989,Han:1988,Tibshirani:2017}.

Suppose further that $K_j=M_j$ is a closed subspace for each $j$ and that $\sum_{j=1}^J M_j^\perp$ is closed.
Then $s_j^{(n)}\in M_j^\perp$ and
\begin{equation*}
\operatorname{proj}_{M_j}(v+s_j^{(n)})
=
\operatorname{proj}_{M_j}v.
\end{equation*}
Hence, the Dykstra iteration in \cref{Alg:Dykstra} reduces to the von Neumann iteration in \cref{Alg:MSC_von_Neumann}.
Moreover, since $\chi_{M_j}^*=\chi_{M_j^\perp}$, the product dual~\eqref{POCS_subsets_dual} reduces to
\begin{equation*}
\min_{q_j\in M_j^\perp}
\frac12\left\|\sum_{j=1}^Jq_j-f\right\|^2.
\end{equation*}
This is the expanded system associated with the linear problem
\begin{equation*}
\min_{q\in M^\perp}\frac12\|q-f\|^2,
\quad \text{ where } \quad
M=\bigcap_{j=1}^JM_j,\
M^\perp=\sum_{j=1}^JM_j^\perp,
\end{equation*}
which is equivalent to the variational problem~\eqref{von_Neumann_variational}.
Therefore, the equivalence between block Gauss--Seidel and SSC in \cref{Sec:MSC} shows that the specialization above recovers the SSC and von Neumann relation established in~\cite{XZ:2002}.

\begin{remark}
\label{Rem:Kaczmarz}
The Kaczmarz algorithm is an iterative projection method for solving linear systems~\cite{LK:2016,SV:2009}.
For a consistent system $Av=f$, where $A\in\mathbb{R}^{m\times d}$ has nonzero rows $a_i$, each equation defines the affine hyperplane
$K_i=\{v\in\mathbb R^d:(a_i,v)=f_i\}$.
Successively projecting onto these hyperplanes gives the Kaczmarz algorithm in \cref{Alg:Kaczmarz}.
This successive projection procedure is a special case of the von Neumann algorithm in \cref{Alg:MSC_von_Neumann}, but with affine subspaces.
The Kaczmarz algorithm is widely used in signal and image reconstruction, particularly in computed tomography~\cite{GBH:1970}.
\end{remark}

\begin{algorithm}
\caption{Kaczmarz algorithm for $Av=f$}
\label{Alg:Kaczmarz}
\begin{algorithmic}[]
\STATE Choose $u^{(0)}\in\mathbb{R}^d$.
\FOR{$n=0,1,2,\dots$}
    \FOR{$i=1,2,\dots,m$}
        \STATE \vspace{-0.5cm}
        \begin{equation*}
        u^{(n+\frac{i}{m})}
        =
        u^{(n+\frac{i-1}{m})}
        +
        \frac{f_i-(a_i,u^{(n+\frac{i-1}{m})})}{\|a_i\|^2}a_i
        \end{equation*}
        \vspace{-0.3cm}
    \ENDFOR
\ENDFOR
\end{algorithmic}
\end{algorithm}

\begin{remark}
\label{Rem:acceleration_GS}
Acceleration can be introduced through the corresponding dual subspace correction methods.
For affine subspaces, PSC may be used as a preconditioner for the conjugate gradient method, including the Bramble--Pasciak--Xu preconditioner~\cite{BPX:1990} and overlapping Schwarz preconditioners~\cite{Park:2024b,TW:2005}; convergence results are given in~\cite{GV:2013,Saad:2003}.
Nesterov-type acceleration of nonlinear PSC was considered in~\cite{LP:2019,PX:2023}.
Acceleration of two-block Gauss--Seidel methods was studied in~\cite{CP:2015,ST:2016} and applied to Dykstra in~\cite{CTV:2017}.
An accelerated randomized Dykstra algorithm was also considered in~\cite{NF:2022}.
\end{remark}

\subsection{Applications}
We next present several applications of the dualization relation between subspace correction and operator splitting methods.

\subsubsection{Rudin--Osher--Fatemi model}
\label{Subsubsec:ROF}
The Rudin--Osher--Fatemi model, originally introduced in~\cite{ROF:1992}, is a fundamental model in mathematical imaging and is widely used for signal and image denoising; see, e.g.,~\cite{CP:2016a,LP:2020}.
For simplicity, we consider the discrete one-dimensional case.
Let $V=\mathbb{R}^d$ be the space of one-dimensional signals, and define the discrete gradient operator $D\colon V\to V$ by
\begin{equation*}
(Dv)_i
=
\begin{cases}
v_{i+1}-v_i,
& \text{if }1\leq i\leq d-1,\\
0,
& \text{if }i=d.
\end{cases}
\end{equation*}
Given a noisy signal $f\in V$ and a regularization parameter
$\alpha>0$, the discrete Rudin--Osher--Fatemi problem is given by
\begin{equation}
\label{ROF_primal}
\min_{v\in V}
\left\{
\frac{\alpha}{2}\|v-f\|^2+\|Dv\|_1
\right\}.
\end{equation}
Problem~\eqref{ROF_primal} is a special case of~\eqref{multiple_convex_primal}, and hence of~\eqref{model_problem_2}. Its dual problem is given by~\cite{Chambolle:2004,KH:2004}
\begin{equation}
\label{ROF_dual}
\min_{q\in V^*}
\left\{
\frac{1}{2\alpha}\|D^*q\|^2-(q,Df)
\right\}
\quad\text{subject to}\quad
\|q\|_\infty\leq1.
\end{equation}
Problem~\eqref{ROF_dual} is also an instance of~\eqref{model_problem_2}, with the constraint represented by an indicator function.
Although~\eqref{ROF_primal} is strongly convex, Schwarz-type methods applied directly to this primal problem may fail to converge because of the nonsmooth term~\cite[Example~6.1]{LN:2017}; see also~\cite[Section~3]{WK:2003}.
The dual problem is not strongly convex because $D^*$ has a nontrivial null space, but its coordinatewise constraint is well suited to domain decomposition methods~\cite{CTWY:2015,HL:2015,LPP:2019,LP:2019}.
Dualizing such methods yields convergent domain decomposition methods for the primal problem~\cite{LP:2020}.

In the dual problem~\eqref{ROF_dual}, dual functionals are represented by their coefficient vectors with respect to the standard Euclidean basis.
Let
\begin{equation*}
\Omega_1 = \{ 1, \dots, d_1 \}, \quad
\Omega_2 = \{ d_1 + 1, \dots, d \}, \quad
\Gamma = \{ d_1 + 1 \},
\end{equation*}
for some $1 \leq d_1 < d$.
With this representation, the dual variable space admits the nonoverlapping decomposition
\begin{equation}
\label{ROF_space_decomposition}
V^* = V_1^* + V_2^*,
\quad V_j^* = \{ q \in V^* : \operatorname{supp} q \subset \Omega_j \}, \quad j = 1,2.
\end{equation}
Subspace correction methods~(\cref{Alg:PSC,Alg:SSC}) for~\eqref{ROF_dual} based on the dual-space decomposition~\eqref{ROF_space_decomposition} were considered and analyzed in~\cite{HL:2015,LP:2019,Park:2020}.

By \cref{Thm:multiple_convex_PR,Thm:multiple_convex_parallel}, \cref{Alg:multiple_convex_PR,Alg:multiple_convex_parallel} applied to~\eqref{ROF_primal} are dualizations of SSC and PSC applied to~\eqref{ROF_dual}, respectively.
For the parallel Douglas--Rachford algorithm, the
$\hat{u}_j^{(n+1)}$-subproblem~\eqref{Alg1:multiple_convex_parallel} is given by
\begin{equation}
\label{ROF_primal_DD_1}
\hat{u}_j^{(n+1)}
\in
\operatornamewithlimits{\arg\min}_{v\in V}
\left\{
\frac{\alpha}{2}\|v-u^{(n)}\|^2
+
(r_j^{(n)},v)
+
\|(Dv)|_{\Omega_j}\|_1
\right\}.
\end{equation}
In~\eqref{ROF_primal_DD_1}, only the degrees of freedom in the subspace $W_j$ are involved in the computation of $(Dv)|_{\Omega_j}$, where
\begin{equation*}
W_1 = \{ w \in V : \operatorname{supp} w \subset \Omega_1 \cup \Gamma \},
\quad W_2 = \{ w \in V : \operatorname{supp} w \subset \Omega_2 \}.
\end{equation*}
Under the dualization relation, $r_j^{(n)}$ is also supported in $W_j$.
Consequently, the minimizer has the form
\[
\hat{u}_j^{(n+1)}
=
u^{(n)}+w_j^{(n+1)},
\quad
w_j^{(n+1)}\in W_j,
\]
and~\eqref{ROF_primal_DD_1} reduces to a local problem on $W_j$.
The Peaceman--Rachford subproblem has the same structure, with $u^{(n+\frac{j-1}{J})}$ replacing $u^{(n)}$.
Hence, both splitting algorithms are domain decomposition methods based on the overlapping decomposition $V=W_1+W_2$.
A direct calculation shows that they are equivalent to the methods proposed in~\cite{LG:2019,LN:2017}, providing a concrete instance of the dualization relation between operator splitting and subspace correction.

\subsubsection{Multinomial logistic regression}
\label{Subsubsec:logistic_regression}
Logistic regression is a fundamental classification model and commonly appears as the final softmax layer of a neural network; see, e.g.,~\cite{HX:2019,HZRS:2016,SZ:2015}.
We consider multinomial logistic regression for $k$-class classification.
Let
\[
\{(x_j,y_j)\}_{j=1}^N
\subset
\mathbb{R}^d\times\{1,\dots,k\}
\]
be a labeled dataset, and let
$w_i\in\mathbb{R}^d$ and $b_i\in\mathbb{R}$ denote the weight and bias associated with class $i$.
Let $V=(\mathbb{R}^d\times\mathbb{R})^k$ and let $\alpha>0$ be a regularization parameter.
The regularized multinomial logistic-regression problem is given by
\begin{equation}
\label{logistic_regression_primal}
\min_{\theta=((w_i,b_i))_{i=1}^k\in V}
\left\{
\frac{1}{N}
\sum_{j=1}^N
\left[
\log
\left(
\sum_{i=1}^k
e^{(w_i,x_j)+b_i}
\right)
-
((w_{y_j},x_j)+b_{y_j})
\right]
+
\frac{\alpha}{2}\|\theta\|^2
\right\}.
\end{equation}
Represent $\theta$ by $\operatorname{col}(w_1,b_1,\dots,w_k,b_k)\in\mathbb{R}^{(d+1)k}$, and define
\begin{equation*}
X_j
=
I_k\otimes
\begin{bmatrix}
x_j\\
1
\end{bmatrix}
\in
\mathbb{R}^{(d+1)k\times k},
\end{equation*}
and $\hat{x}\in\mathbb{R}^{(d+1)k}$ by
\begin{equation*}
(\hat{x},\theta)
=
\sum_{j=1}^N
((w_{y_j},x_j)+b_{y_j})
=
\sum_{i=1}^k
\sum_{j:y_j=i}
((w_i,x_j)+b_i).
\end{equation*}
Since $X_j^*\theta=((w_i,x_j)+b_i)_{i=1}^k$, multiplying the objective in~\eqref{logistic_regression_primal} by $N$ shows that it is an instance of~\eqref{multiple_convex_primal}, and hence of~\eqref{model_problem_2}, with
\begin{equation*}
F(\theta)
=
\frac{N\alpha}{2}\|\theta\|^2-(\hat{x},\theta),
\quad
G_j=\operatorname{LSE}_k,
\quad
B_j=X_j^*,
\quad 1\leq j\leq N.
\end{equation*}
Here, $\operatorname{LSE}_k$ is defined in \cref{Table:conjugate}.
The corresponding dual problem is given by
\begin{equation}
\label{logistic_regression_dual}
\min_{q_j\in\mathbb{R}^k}
\left\{
\frac{1}{2N\alpha}
\left\|
\sum_{j=1}^N X_jq_j-\hat{x}
\right\|^2
+
\sum_{j=1}^N
\operatorname{LSE}_k^*(q_j)
\right\}.
\end{equation}
Problem~\eqref{logistic_regression_dual} is also an instance of~\eqref{model_problem_2}, with $\undertilde{q}=(q_1,\dots,q_N)\in\mathbb{R}^{Nk}$.
If $\theta$ and $\undertilde{p}=(p_1,\dots,p_N)\in\mathbb{R}^{Nk}$ are primal and dual solutions, respectively, then the primal--dual relations are given by
\begin{equation*}
-\sum_{j=1}^N X_jp_j
=
N\alpha\theta-\hat{x},
\quad
p_j
=
(\operatorname{LSE}_k)'(X_j^*\theta),
\quad
1\leq j\leq N.
\end{equation*}

Each block $q_j$ in~\eqref{logistic_regression_dual} corresponds to one data pair $(x_j,y_j)$, so updating a single block does not require recomputing the entire gradient.
This structure is well suited to block coordinate descent methods~\cite{FR:2015,FR:2016,Nesterov:2012,Wright:2015}, including those studied in~\cite{RT:2016,SZ:2013,SZ:2014,YHL:2011}.

Taking $J=N$ and $B_j=X_j^*$ in \cref{Alg:multiple_convex_PR}, and writing $\phi_j^{(n)}\in\mathbb{R}^{(d+1)k}$ for the Riesz representative of $r_j^{(n)}\in (\mathbb{R}^{(d+1)k})^*$, gives the following iteration for~\eqref{logistic_regression_primal}, which processes one data pair at a time:
\begin{subequations}
\label{Alg:logistic_regression_PR}
\begin{align}
\label{Alg1:logistic_regression_PR}
\theta^{(n+\frac{j}{N})}
&=
\operatornamewithlimits{\arg\min}_{\theta \in\mathbb{R}^{(d+1)k}}
\left\{
\frac{N \alpha}{2}
\left\|
\theta-\theta^{(n+\frac{j-1}{N})}
\right\|^2
+
( \phi_j^{(n)},\theta)
+
\operatorname{LSE}_k(X_j^*\theta)
\right\}, \\
\phi_j^{(n+1)}
&=
\phi_j^{(n)}
+
N\alpha
\left(
\theta^{(n+\frac{j}{N})}
-
\theta^{(n+\frac{j-1}{N})}
\right),
\quad 1\leq j\leq N.
\end{align}
\end{subequations}
Completing the square in~\eqref{Alg1:logistic_regression_PR} shows that this is a correction-based incremental proximal iteration.
Moreover, by \cref{Thm:multiple_convex_PR}, iteration~\eqref{Alg:logistic_regression_PR} is a dualization of the dual coordinate descent method proposed in~\cite{YHL:2011}, which is SSC applied to~\eqref{logistic_regression_dual}.

The loop over $j$ can instead be parallelized by applying \cref{Alg:multiple_convex_parallel} to~\eqref{logistic_regression_primal}, or randomized by selecting $j$ at each iteration.
The randomized method is a dualization of the stochastic dual coordinate ascent algorithm proposed in~\cite{SZ:2013,SZ:2014}; see also~\cite{GO:2012,HXZ:2019,JPX:2026} for randomized subspace correction methods.
Parallelization and randomization can be combined as in~\cite{FR:2015,FR:2016,RT:2016}, yielding a stochastic mini-batch scheme closely related to the mini-batch principle used in stochastic gradient descent and other stochastic optimization methods for training deep learning models~\cite{GBC:2016}.
  
\subsubsection{Federated learning}
Federated learning is a distributed training paradigm in which data are stored on a collection of clients, while a central server maintains a global model parameter.  The clients perform local computations using their own datasets and communicate model updates or local solutions to the server, rather than sharing raw data.  This setting has become important as data are increasingly decentralized and privacy considerations limit direct data sharing~\cite{Kairouz:2021,LSTS:2020,MMRHA:2017}.  Since communication between clients and the server is often expensive, a central issue in federated learning is to design algorithms that use local training effectively while controlling the effect of data heterogeneity.

Let $\theta\in V$ denote the global model parameter, and let $G_j\colon V\to\overline{\mathbb R}$ be the local objective associated with the dataset $\mathcal D_j$ on client $j$. With $v=\theta$, the regularized federated learning model is precisely~\eqref{multiple_convex_special}; factors such as $1/J$ may be absorbed into $G_j$. Its dual is the specialization of~\eqref{multiple_convex_dual} with $F(\theta)=(\alpha/2)\|\theta\|^2$ and $B_j=I$:
\begin{equation}
\label{FL_application_dual}
\min_{q_j\in V^*}
\left\{
\sum_{j=1}^J G_j^*(q_j)
+
\frac{1}{2\alpha}
\left\|
\sum_{j=1}^J q_j
\right\|^2
\right\}.
\end{equation}
Problem~\eqref{FL_application_dual} is also an instance of~\eqref{model_problem_2}, and the primal--dual relation~\eqref{multiple_primal_dual_relations} reduces to
\begin{equation*}
\theta=-\frac{1}{\alpha}\sum_{j=1}^J p_{j},
\quad
p_{j}\in\partial G_j(\theta),
\quad 1\leq j\leq J.
\end{equation*}
Thus the server recovers $\theta$ from the client-wise dual components $p_j$. This relation is analogous to that in dual coordinate methods for regularized empirical risk minimization~\cite{SZ:2013,SZ:2014}; one may also refer to the duality-based federated learning approach in~\cite{PX:2023}.

Applying the parallel splitting algorithm in~\cref{Alg:multiple_convex_parallel} to~\eqref{multiple_convex_special}, together with the update in~\cref{Rem:multiple_convex_parallel}, gives the independent client subproblems
\begin{equation}
\label{FL_application_local}
\hat\theta_j^{(n+1)}
=
\operatornamewithlimits{\arg\min}_{\theta\in V}
\left\{
\frac{\alpha}{2}\|\theta-\theta^{(n)}\|^2
+
(r_j^{(n)},\theta)
+
G_j(\theta)
\right\},
\quad 1\leq j\leq J.
\end{equation}
Each subproblem uses only $G_j$ and can be solved on client $j$. The server then applies the aggregation formula in~\cref{Rem:multiple_convex_parallel}:
\begin{equation*}
\theta^{(n+1)}=(1-\tau J)\theta^{(n)}+\tau\sum_{j=1}^J\hat\theta_j^{(n+1)}.
\end{equation*}
For $\tau=1/J$, this is the average of the client updates. For heterogeneous data, client drift is often mitigated by control-variate or gradient-shift corrections~\cite{KKMRSS:2020,KMR:2020,MMSR:2022}. The linear term $(r_j^{(n)},\theta)$ in~\eqref{FL_application_local} resembles such corrections, but here it follows from dualization: by~\cref{Thm:multiple_convex_parallel}, $r_j^{(n)}=-p_j^{(n)}$ because $B_j=I$. It is therefore dictated by the primal--dual consistency between~\eqref{FL_application_dual} and~\eqref{multiple_convex_special}.
  
\section{Alternating direction methods of multipliers (ADMMs)}
\label{Sec:ADMM}

In this section, we study standard ADMM and its variants for solving multi-block constrained optimization problems.
We also revisit the classical connection between two-block ADMM and Douglas--Rachford splitting within the dualization framework and use it to derive a new multi-block ADMM-type method via dualization.

\subsection{Plain method and its variants}
\label{Subsec:ADMM_plain}

As a multi-block case of~\eqref{model_problem_3}, we consider the linearly constrained convex optimization problem
\begin{equation}
\label{ADMM_multiple}
\min_{v_j\in V_j}
\sum_{j=1}^J F_j(v_j)
\quad\text{subject to}\quad
\sum_{j=1}^J B_jv_j=g,
\end{equation}
where each $V_j$, $1\leq j\leq J$, and $W$ are real Hilbert spaces, and $g \in W$.
Each function $F_j\colon V_j\to\overline{\mathbb{R}}$ is proper, convex,
and lower semicontinuous, and each operator $B_j\colon V_j\to W$ is bounded and linear.
In this section, we write
\begin{equation*}
\undertilde{V}=\prod_{j=1}^J V_j,
\quad
\undertilde{v}=(v_1,\dots,v_J)\in\undertilde{V}.
\end{equation*}
The same convention is used for iterates, e.g.,
$\undertilde{u}^{(n)}=(u_1^{(n)},\dots,u_J^{(n)})\in\undertilde{V}$.

For a block-structured primal variable, ADMM modifies the augmented Lagrangian method~\eqref{Alg:ALM} by replacing the minimization with respect to the full primal variable in~\eqref{Alg1:ALM} with blockwise minimization, while retaining the multiplier update~\eqref{Alg2:ALM}.

Let $\mathcal{L}_{\beta}(\undertilde{v},\lambda)$ denote the augmented Lagrangian for~\eqref{ADMM_multiple} (cf.~\eqref{ALM_augmented_Lagrangian}):
\begin{equation*}
\mathcal{L}_{\beta} (\undertilde{v}, \lambda)
=
\sum_{j=1}^J F_j ( v_j )
+
\left( \lambda, \sum_{j=1}^J B_j v_j - g \right)
+
\frac{\beta}{2}
\left\|
\sum_{j=1}^J B_j v_j - g
\right\|^2.
\end{equation*}
The plain ADMM~(see, e.g.,~\cite{CHYY:2014,WYZ:2019}) for solving~\eqref{ADMM_multiple} is presented in \cref{Alg:ADMM_multiple}.

\begin{algorithm}
\caption{ADMM for~\eqref{ADMM_multiple}~\cite{CHYY:2014,WYZ:2019}}
\label{Alg:ADMM_multiple}
\begin{algorithmic}[]
\STATE Given $\beta > 0$:
\STATE Choose $\undertilde{u}^{(0)} \in \undertilde{V}$ and $\lambda^{(0)} \in W^*$.
\FOR{$n=0,1,2,\dots$}
\FOR{$j=1, 2, \dots, J$}
\STATE \vspace{-0.5cm} \begin{equation*}
    u_j^{(n+1)}
    \in
    \operatornamewithlimits{\arg\min}_{v_j \in V_j}
    \mathcal{L}_{\beta}
    (u_1^{(n+1)}, \dots, u_{j-1}^{(n+1)}, v_j,
    u_{j+1}^{(n)}, \dots, u_J^{(n)}, \lambda^{(n)})
    \end{equation*}
    \vspace{-0.3cm}
\ENDFOR
\STATE \vspace{-0.5cm} \begin{equation*}
    \lambda^{(n+1)}
    =
    \lambda^{(n)}
    +
    \beta
    \left(
    \sum_{j=1}^{J} B_j u_j^{(n+1)} - g
    \right)
    \end{equation*}
    \vspace{-0.3cm}
\ENDFOR
\end{algorithmic}
\end{algorithm}

Under the standard solvability assumptions in~\cite{HY:2014,NLRPJ:2015}, the convergence of \cref{Alg:ADMM_multiple} is guaranteed in the two-block case, i.e., when $J = 2$. In contrast, it is not guaranteed in general for $J \geq 3$, as shown in~\cite{CHYY:2014}.
To ensure convergence of \cref{Alg:ADMM_multiple} for $J \geq 3$, additional assumptions are required; relevant conditions are discussed in~\cite{TY:2018} and the references therein.

Other convergent modifications use parallel multi-block or relaxed full-Jacobian updates~\cite{DLPY:2017,HHY:2015}.
One remedy for the nonconvergence of plain multi-block ADMM is symmetrization~\cite{CST:2017,XXY:2017}.
The symmetrized ADMM is presented in \cref{Alg:ADMM_multiple_symmetrized}; its convergence has been studied in the context of quadratic optimization in~\cite{CST:2017,XXY:2017}.
Unlike the plain method, the update step for the primal variables consists of two sweeps: one from $1$ to $J$, followed by another from $J$ to $1$.

\begin{algorithm}
\caption{Symmetrized ADMM for~\eqref{ADMM_multiple}~\cite{CST:2017,XXY:2017}}
\label{Alg:ADMM_multiple_symmetrized}
\begin{algorithmic}[]
\STATE Given $\beta > 0$:
\STATE Choose $\undertilde{u}^{(0)} \in \undertilde{V}$ and $\lambda^{(0)} \in W^*$.
\FOR{$n=0,1,2,\dots$}
\FOR{$j=1, 2, \dots, J$}
\STATE \vspace{-0.5cm} \begin{equation*}
    u_j^{(n + \frac{1}{2})}
    \in
    \operatornamewithlimits{\arg\min}_{v_j \in V_j}
    \mathcal{L}_{\beta}
    (u_1^{(n+\frac{1}{2})}, \dots, u_{j-1}^{(n+\frac{1}{2})},
    v_j, u_{j+1}^{(n)}, \dots, u_J^{(n)}, \lambda^{(n)})
    \end{equation*}
    \vspace{-0.3cm}
\ENDFOR
\FOR{$j=J, J-1, \dots, 1$}
\STATE \vspace{-0.5cm} \begin{equation*}
    u_j^{(n + 1)}
    \in
    \operatornamewithlimits{\arg\min}_{v_j \in V_j}
    \mathcal{L}_{\beta}
    (u_1^{(n+\frac{1}{2})}, \dots, u_{j-1}^{(n+\frac{1}{2})},
    v_j, u_{j+1}^{(n+1)}, \dots, u_J^{(n+1)}, \lambda^{(n)})
    \end{equation*}
    \vspace{-0.3cm}
\ENDFOR
\STATE \vspace{-0.5cm} \begin{equation*}
        \lambda^{(n+1)}
        =
        \lambda^{(n)}
        +
        \beta
        \left(
        \sum_{j=1}^{J} B_j u_j^{(n+1)} - g
        \right)
    \end{equation*}
    \vspace{-0.3cm}
\ENDFOR
\end{algorithmic}
\end{algorithm}

We now interpret the plain and symmetrized ADMMs in the context of quadratic optimization.
Specifically, we show that, in the quadratic setting, these methods correspond to instances of the inexact Uzawa method introduced in~\cite{BPV:1997}, with different primal smoothers.

We consider the following constrained quadratic optimization problem, which is a special case of~\eqref{ADMM_multiple}, and hence an instance of~\eqref{model_problem_3}:
\begin{equation}
\label{ADMM_linear}
\min_{v_j \in V_j}
\sum_{j=1}^J
\left[
\frac{1}{2} (A_j v_j,v_j)
-
(f_j,v_j)
\right]
\quad \text{ subject to } \quad
\sum_{j=1}^J B_j v_j = g,
\end{equation}
where each $A_j \colon V_j \to V_j^*$ is the linear map induced by an SPD bilinear form on $V_j$ and $f_j \in V_j^*$.
The optimality conditions for the corresponding augmented Lagrangian formulation form the saddle point system
\begin{equation}
\label{ADMM_linear_saddle}
\begin{bmatrix}
\undertilde{A}_{\beta} & \undertilde{B}^* \\
\undertilde{B} & 0
\end{bmatrix}
\begin{bmatrix}
\undertilde{v} \\
\lambda
\end{bmatrix}
=
\begin{bmatrix}
\undertilde{f}_{\beta} \\
g
\end{bmatrix},
\end{equation}
where
\begin{equation*}
\undertilde{A}_{\beta}
=
\undertilde{A}
+
\beta \undertilde{B}^*\undertilde{B},
\quad
\undertilde{A}
=
\operatorname{diag} (A_j),
\quad
\undertilde{B}
=
[B_1, \dots, B_J],
\quad
\undertilde{f}_{\beta}
=
\begin{bmatrix}
f_1+\beta B_1^* g \\
\vdots \\
f_J+\beta B_J^* g
\end{bmatrix}.
\end{equation*}
Let $\undertilde{D}_{\beta}$ and $\undertilde{L}_{\beta}$ denote the block diagonal and strictly block lower-triangular parts of $\undertilde{A}_{\beta}$, respectively, so that
\begin{equation*}
\undertilde{A}_{\beta}
=
\undertilde{L}_{\beta}
+
\undertilde{D}_{\beta}
+
\undertilde{L}_{\beta}^*.
\end{equation*}

In~\cite[Algorithm~2.3]{BPV:1997}, a general class of iterative methods for solving saddle point systems of the form~\eqref{ADMM_linear_saddle}, known as the inexact Uzawa method, was introduced:
\begin{subequations}
\label{ADMM_linear_inexact_Uzawa}
\begin{align}
\undertilde{u}^{\mathrm{new}}
&=
\undertilde{u}^{\mathrm{old}}
+
R_{\undertilde{V}}
(
\undertilde{f}_{\beta}
-
(
\undertilde{A}_{\beta} \undertilde{u}^{\mathrm{old}}
+
\undertilde{B}^* \lambda^{\mathrm{old}}
)
),
\\
\lambda^{\mathrm{new}}
&=
\lambda^{\mathrm{old}}
+
R_W
(
\undertilde{B} \undertilde{u}^{\mathrm{new}} - g
),
\end{align}
\end{subequations}
where $R_{\undertilde{V}} \colon \undertilde{V}^* \to \undertilde{V}$ and $R_W \colon W \to W^*$ are boundedly invertible linear operators.

For the plain ADMM~(\cref{Alg:ADMM_multiple}), alternating minimization with respect to $\undertilde{v}$ is equivalent to block Gauss--Seidel smoothing for $\undertilde{A}_{\beta}$~(cf.~\cite{Saad:2003,XZ:2017}), and hence
\begin{equation}
\label{Uzawa_ADMM_multiple}
R_{\undertilde{V}}
=
(\undertilde{L}_{\beta}+\undertilde{D}_{\beta})^{-1},
\quad
R_W
=
\beta I.
\end{equation}
For the symmetrized ADMM~(\cref{Alg:ADMM_multiple_symmetrized}), the two sweeps correspond to symmetrized block Gauss--Seidel smoothing for $\undertilde{A}_{\beta}$, giving
\begin{equation*}
R_{\undertilde{V}}
=
((\undertilde{L}_{\beta}+\undertilde{D}_{\beta})^*)^{-1}
\undertilde{D}_{\beta}
(\undertilde{L}_{\beta}+\undertilde{D}_{\beta})^{-1},
\quad
R_W
=
\beta I.
\end{equation*}

A sufficient condition for the convergence of~\eqref{ADMM_linear_inexact_Uzawa}, given in~\cite[Corollary~3.2]{BPV:1997}, states that $R_{\undertilde{V}}^{-1}\colon \undertilde{V}\to \undertilde{V}^*$ and $R_W\colon W\to W^*$ are SPD linear maps and satisfy
\begin{subequations}
\begin{align}
\label{ADMM_linear_inexact_Uzawa_sufficient1}
( \undertilde{A}_{\beta} \undertilde{v}, \undertilde{v} )
\leq
( R_{\undertilde{V}}^{-1} \undertilde{v}, \undertilde{v} ),
&\quad
\undertilde{v} \in \undertilde{V},
\\
\label{ADMM_linear_inexact_Uzawa_sufficient2}
( \lambda, \undertilde{B} \undertilde{A}_{\beta}^{-1} \undertilde{B}^* \lambda )
\leq
( \lambda, R_W^{-1} \lambda ),
&\quad
\lambda \in W^*.
\end{align}
\end{subequations}
The symmetrized ADMM satisfies~\eqref{ADMM_linear_inexact_Uzawa_sufficient1} by the standard argument for symmetrized block Gauss--Seidel smoothing applied to SPD systems.
Moreover, since $R_W=\beta I$, condition~\eqref{ADMM_linear_inexact_Uzawa_sufficient2} becomes
\begin{equation*}
( \lambda, \undertilde{B} \undertilde{A}_{\beta}^{-1} \undertilde{B}^* \lambda )
\leq
\beta^{-1}
( \lambda, \lambda ),
\quad
\lambda \in W^*,
\end{equation*}
which follows from the corresponding operator inequality.
Therefore, the convergence of the symmetrized ADMM for~\eqref{ADMM_linear} follows from the inexact Uzawa theory.

In contrast, for the plain ADMM in~\eqref{Uzawa_ADMM_multiple}, the operator $R_{\undertilde{V}}$ is not symmetric, and the above sufficient condition does not apply.
This explains, from the inexact Uzawa viewpoint, why the convergence of the plain multi-block ADMM cannot be ensured in general.

\begin{remark}
\label{Rem:randomized_ADMM}
Another approach to the nonconvergence issue is randomization, where the order of the block updates is randomized; see, e.g.,~\cite{MZY:2020,SLY:2015,SLY:2020}.
Connections between randomized block orderings and averaged Gauss--Seidel-type smoothers have also been considered in~\cite{CG:2023}.
\end{remark}

\subsection{Two-block ADMM and Douglas--Rachford splitting}
In the two-block case~($J = 2$), the equivalence between ADMM and the Douglas--Rachford splitting algorithm applied to a dual problem is classical; see, e.g.,~\cite{EB:1992,Setzer:2011}.
Our purpose in this subsection is therefore not to present this equivalence as new, but to express it in the notation of \cref{Def:dualization}.
This formulation will then be used as a template for the multi-block constructions later.

We assume that the infimal postcomposition appearing below is proper and lower semicontinuous and that the infimum in its definition is attained.

We consider the following two-block constrained optimization problem, which is the $J = 2$ instance of~\eqref{ADMM_multiple}, and hence an instance of~\eqref{model_problem_3}:
\begin{equation}
\label{ADMM_two}
\min_{v_1 \in V_1,\  v_2 \in V_2 } \left\{ F_1 (v_1) + F_2 (v_2) \right\}
\quad \text{ subject to } \quad
B_1 v_1 +  B_2 v_2 = g.
\end{equation}
Each iteration of ADMM~(\cref{Alg:ADMM_multiple}) for solving~\eqref{ADMM_two} is presented as follows:
\begin{subequations}
\label{Alg:ADMM_two}
\begin{align}
        \label{Alg1:ADMM_two}
        &u_1^{(n+1)} \in \operatornamewithlimits{\arg\min}_{v_1 \in V_1} \left\{ F_1 ( v_1 )  + (\lambda^{(n)},B_1v_1) + \frac{\beta}{2} \| B_1 v_1 + B_2 u_2^{(n)} - g \|^2 \right\} \\
        \label{Alg2:ADMM_two}
        &u_2^{(n+1)} \in \operatornamewithlimits{\arg\min}_{v_2 \in V_2} \left\{ F_2 ( v_2 ) + (\lambda^{(n)},B_2v_2) + \frac{\beta}{2} \| B_1 u_1^{(n+1)} + B_2 v_2 - g \|^2 \right\} \\
        \label{Alg3:ADMM_two}
        &\lambda^{(n+1)} = \lambda^{(n)} + \beta ( B_1 u_1^{(n+1)} + B_2 u_2^{(n+1)} - g )
    \end{align}
\end{subequations}

We observe that the problem~\eqref{ADMM_two} is equivalent to the following unconstrained convex optimization problem:
\begin{equation}
\label{ADMM_two_primal}
\min_{v_1 \in V_1} \left\{ F_1 ( v_1 ) + ((- B_2) \triangleright F_2) ( B_1 v_1 - g ) \right\},
\end{equation}
where $\triangleright$ denotes the infimal postcomposition defined in~\eqref{infimal_postcomposition}.
Problem~\eqref{ADMM_two_primal} is an instance of~\eqref{model_problem_2}.
Applying Fenchel--Rockafellar duality, we obtain the following dual problem~(cf.~\eqref{infimal_postcomposition_conjugate}):
\begin{equation}
\label{ADMM_two_dual}
    \min_{q \in W^*} \left\{ F_1^* (- B_1^* q) + F_2^* (- B_2^* q) + (q,g) \right\}.
\end{equation}
Problem~\eqref{ADMM_two_dual} is also an instance of~\eqref{model_problem_2}.
The first primal--dual relation~\eqref{primal_dual_relation_1} is written as
\begin{equation*}
- B_1^* p \in \partial F_1(u_{1}).
\end{equation*}
The Douglas--Rachford splitting algorithm (see~\eqref{DR_conventional}) for solving the dual problem~\eqref{ADMM_two_dual} is presented as follows:
\begin{subequations}
\label{Alg:ADMM_two_dual}
\begin{align}
    \label{Alg1:ADMM_two_dual}
    &p^{(n+1)} = \operatornamewithlimits{\arg\min}_{q \in W^*} \left\{ F_1^* (- B_1^* q) + (q,g) + \frac{1}{2\beta} \|q+r^{(n)}-2s^{(n)}\|^2 \right\} \\
    \label{Alg2:ADMM_two_dual}
    &r^{(n+1)} = p^{(n+1)} + r^{(n)} - s^{(n)} \\
    \label{Alg3:ADMM_two_dual}
    &s^{(n+1)} = \operatornamewithlimits{\arg\min}_{s \in W^*} \left\{ F_2^* (- B_2^* s) + \frac{1}{2 \beta} \| s - r^{(n+1)} \|^2 \right\}
    \end{align}
\end{subequations}
Here the initial values satisfy $r^{(0)},s^{(0)}\in W^*$.

We next summarize the dualization relation between~\eqref{Alg:ADMM_two} and~\eqref{Alg:ADMM_two_dual}.

\begin{theorem}
\label{Thm:ADMM_two}
For compatible initial data, ADMM~\eqref{Alg:ADMM_two} for solving the primal problem~\eqref{ADMM_two} is a dualization of the Douglas--Rachford splitting algorithm~\eqref{Alg:ADMM_two_dual} for solving the dual problem~\eqref{ADMM_two_dual}.
More precisely, if
\begin{equation*}
- B_2 u_2^{(0)} = \frac{1}{\beta} (r^{(0)} - s^{(0)}), \quad
\lambda^{(0)} = s^{(0)},
\end{equation*}
then we have
\begin{equation*}
- B_1^* p^{(n)} \in \partial F_1 (u_1^{(n)}), \quad
- B_2 u_2^{(n)} = \frac{1}{\beta} (r^{(n)} - s^{(n)}), \quad
\lambda^{(n)} = s^{(n)},
\quad n \geq 1.
\end{equation*}
\end{theorem}
\begin{proof}
Assume that, for some $n\geq 0$,
\begin{equation}
\label{ADMM_two_dual1}
- B_2 u_2^{(n)}
=
\frac{1}{\beta}(r^{(n)}-s^{(n)}),
\quad
\lambda^{(n)}=s^{(n)}.
\end{equation}
Then the $p^{(n+1)}$-subproblem in~\eqref{Alg1:ADMM_two_dual}, after absorbing the term $(q,g)$ into the quadratic term, is Fenchel--Rockafellar dual to the $u_1^{(n+1)}$-subproblem~\eqref{Alg1:ADMM_two}.
The corresponding primal--dual relations give
\begin{subequations}
\begin{align}
\label{ADMM_two_dual2}
B_1 u_1^{(n+1)} + B_2 u_2^{(n)} - g
&=
\frac{1}{\beta}(p^{(n+1)}-\lambda^{(n)}), \\
- B_1^* p^{(n+1)}
&\in
\partial F_1(u_1^{(n+1)}).
\end{align}
\end{subequations}
Using~\eqref{Alg2:ADMM_two_dual}, \eqref{ADMM_two_dual1}, and~\eqref{ADMM_two_dual2}, we obtain
\begin{equation}
\label{ADMM_two_r_update}
r^{(n+1)}
=
\lambda^{(n)}
+
\beta(B_1u_1^{(n+1)}-g).
\end{equation}

With~\eqref{ADMM_two_r_update}, the $s^{(n+1)}$-subproblem in~\eqref{Alg3:ADMM_two_dual} is Fenchel--Rockafellar dual to the $u_2^{(n+1)}$-subproblem~\eqref{Alg2:ADMM_two}.
The first primal--dual relation gives
\begin{equation}
\label{ADMM_two_dual4}
B_1u_1^{(n+1)}+B_2u_2^{(n+1)}-g
=
\frac{1}{\beta}(s^{(n+1)}-\lambda^{(n)}).
\end{equation}
Hence, by~\eqref{Alg3:ADMM_two},
\begin{equation*}
\lambda^{(n+1)}
=
\lambda^{(n)}
+
\beta(B_1u_1^{(n+1)}+B_2u_2^{(n+1)}-g)
=
s^{(n+1)}.
\end{equation*}
Moreover, combining~\eqref{ADMM_two_dual2}, \eqref{ADMM_two_dual4}, and~\eqref{Alg2:ADMM_two_dual}, we get
\begin{equation*}
- B_2u_2^{(n+1)}
=
\frac{1}{\beta}(r^{(n+1)}-s^{(n+1)}).
\end{equation*}
Thus~\eqref{ADMM_two_dual1} holds with $n$ replaced by $n+1$, and the result follows by induction.
\end{proof}

\subsection{New multi-block method derived from dualization}
\label{Subsec:multi-block_ADMM}
The direct extension of ADMM to three or more blocks may fail to converge in general; see, for example,~\cite{CHYY:2014}. Building on the two-block ADMM--Douglas--Rachford relation in the preceding subsection, we derive a dual formulation of~\eqref{ADMM_multiple}, apply the multi-block Peaceman--Rachford method in~\cref{Sec:Splitting}, and dualize the resulting iteration to obtain an ADMM-type method. Whenever convergence is available for the Peaceman--Rachford method, it transfers to the derived method through the dualization relation established below.

We consider the $J$-block constrained problem~\eqref{ADMM_multiple} with $J\geq 3$. We assume that the infimal postcomposition
$(-B_J)\triangleright F_J$ is proper and lower semicontinuous and that
the infimum in its definition is attained.
Eliminating the last block variable shows that~\eqref{ADMM_multiple}
is equivalent to
\begin{equation}
\label{ADMM_multiple_primal}
\min_{v_j\in V_j,\ 1\leq j\leq J-1}
\left\{
\sum_{j=1}^{J-1}F_j(v_j)
+
((-B_J)\triangleright F_J)
\left(
\sum_{j=1}^{J-1}B_jv_j-g
\right)
\right\}.
\end{equation}
Problem~\eqref{ADMM_multiple_primal} is an instance of~\eqref{model_problem_2} with
\begin{multline*}
V\leftarrow\prod_{j=1}^{J-1}V_j,
\quad
F((v_j)_{j=1}^{J-1})
\leftarrow
\sum_{j=1}^{J-1}F_j(v_j),\\
G(w)
\leftarrow
((-B_J)\triangleright F_J)(w-g),
\quad
B\leftarrow[B_1,\dots,B_{J-1}].
\end{multline*}
Under the corresponding qualification condition, Fenchel--Rockafellar
duality and~\eqref{infimal_postcomposition_conjugate} yield the dual problem
\begin{equation}
\label{ADMM_multiple_dual}
\min_{q\in W^*}
\left\{
(q,g)
+
\sum_{j=1}^J F_j^*(-B_j^*q)
\right\}.
\end{equation}
Problem~\eqref{ADMM_multiple_dual} is an instance of~\eqref{model_problem_2}.

If $(u_1,\dots,u_{J-1})$ and $p$ are primal and dual
solutions of~\eqref{ADMM_multiple_primal} and~\eqref{ADMM_multiple_dual}, respectively, then their primal--dual relations include
\begin{equation*}
-B_j^*p\in\partial F_j(u_j),
\quad 1\leq j\leq J-1,
\end{equation*}
and
\begin{equation*}
\sum_{j=1}^{J-1}B_ju_j-g
\in
\partial(F_J^*\circ(-B_J^*))(p).
\end{equation*}
In particular, if $u_J$ is chosen so that
\[
u_J\in\partial F_J^*(-B_J^*p)
\quad\text{and}\quad
\sum_{j=1}^J B_ju_j=g,
\]
then these relations can be written uniformly as
\[
-B_j^*p\in\partial F_j(u_j),
\quad 1\leq j\leq J.
\]

We assume in addition that the functions $F_{J-1}^* \circ (-B_{J-1}^*)$ and $F_J^* \circ (-B_J^*)$ are strongly convex. 
For instance, this holds if $F_{J-1}$ and $F_J$ are $L$-smooth and $B_{J-1}$ and $B_J$ are surjective.
Choose $\beta>0$ so that the functions $H_j$ defined by
\begin{equation}
\label{H_j}
    H_j(q) = \begin{cases}
    F_j^*(-B_j^*q), & \quad 1 \leq j \leq J-2, \\
    F_j^*(-B_j^*q) - \frac{1}{2\beta} \| q \|^2, & \quad j = J-1, J
    \end{cases}
\end{equation}
are convex.
Then the dual problem~\eqref{ADMM_multiple_dual} can be rewritten as
\begin{equation}
\label{ADMM_multiple_dual_modified}
\min_{q\in W^*}
\left\{
    \frac{1}{\beta}\|q\|^2
    +
    (q,g)
    +
\sum_{j=1}^J H_j(q)
\right\}.
\end{equation}
Problem~\eqref{ADMM_multiple_dual_modified} is an instance of~\eqref{model_problem_2}.

Applying the dualization-based Peaceman--Rachford splitting algorithm
(\cref{Alg:multiple_convex_PR}) to~\eqref{ADMM_multiple_dual_modified} gives the iteration
\begin{subequations}
\label{ADMM_multiple_PR}
\begin{align}
    p^{(n+\frac{j}{J})}
    &=
    \operatornamewithlimits{\arg\min}_{q\in W^*}
    \left\{
        \frac{1}{\beta}
        \|q-p^{(n+\frac{j-1}{J})}\|^2
        +
        (q,r_j^{(n)})
        +
        H_j(q)
    \right\},
    \quad 1\leq j\leq J, \\
    r_j^{(n+1)}
    &=
    r_j^{(n)}
    +
    \frac{2}{\beta}
    (
        p^{(n+\frac{j}{J})}
        -
        p^{(n+\frac{j-1}{J})}
    ),
    \quad 1\leq j\leq J.
\end{align}
\end{subequations}
If the initial values $\undertilde r^{(0)}=(r_1^{(0)},\dots,r_J^{(0)})\in W^J$ satisfy
\[
    \sum_{j=1}^J r_j^{(0)}
    =
    \frac{2}{\beta}p^{(0)}+g,
\]
then~\eqref{ADMM_multiple_PR} preserves this relation:
\[
    \sum_{j=1}^J r_j^{(n)}
    =
    \frac{2}{\beta}p^{(n)}+g,
    \quad n\geq 0.
\]

Dualizing~\eqref{ADMM_multiple_PR} yields the following sequential ADMM-type method.

\begin{algorithm}
\caption{Dualization-based ADMM for~\eqref{ADMM_multiple}}
\label{Alg:ADMM_multiple_primal}
\begin{algorithmic}[]
\STATE Given $\beta>0$:
\STATE Set $\gamma_j=\frac{1}{2}$ for $1\leq j\leq J-2$ and $\gamma_j=1$ for $j=J-1,J$.
\STATE Choose $\undertilde{s}^{(0)}=(s_1^{(0)},\dots,s_J^{(0)})\in W^J$ and $\lambda^{(0)}\in W^*$.
\FOR{$n=0,1,2,\dots$}
    \FOR{$j=1,2,\dots,J$}
        \STATE
        \vspace{-0.5cm}
        \begin{equation*}
        \resizebox{0.9\textwidth}{!}{$\displaystyle
        \begin{aligned}
        \hat{u}_j^{(n+1)}
        &\in
        \operatornamewithlimits{\arg\min}_{v_j\in V_j}
        \left\{
            F_j(v_j)
            +
            2\gamma_j
            (
                \lambda^{(n+\frac{j-1}{J})},
                B_j v_j-s_j^{(n)}
            )
            +
            \frac{\gamma_j\beta}{2}
            \|
                B_j v_j-s_j^{(n)}
            \|^2
        \right\} \\
        \lambda^{(n+\frac{j}{J})}
        &=
        2\gamma_j\lambda^{(n+\frac{j-1}{J})}
        +
        \gamma_j\beta
        (
            B_j\hat{u}_j^{(n+1)}
            -
            s_j^{(n)}
        ) \\
        s_j^{(n+1)}
        &=
        s_j^{(n)}
        +
        \frac{2}{\beta}
        (
            \lambda^{(n+\frac{j}{J})}
            -
            \lambda^{(n+\frac{j-1}{J})}
        )
        \end{aligned}
        $}
        \end{equation*}
        \vspace{-0.3cm}
    \ENDFOR
    \STATE Choose $\undertilde{u}^{(n+1)}=(u_1^{(n+1)},\dots,u_J^{(n+1)})\in\undertilde{V}$ such that
    \STATE \vspace{-0.5cm}
    \begin{equation*}
    -B_j^*\lambda^{(n+1)}\in\partial F_j(u_j^{(n+1)}),
    \quad 1\leq j\leq J.
    \end{equation*}
    \vspace{-0.3cm}
\ENDFOR
\end{algorithmic}
\end{algorithm}

Unlike the direct multi-block ADMM~(\cref{Alg:ADMM_multiple}), 
\Cref{Alg:ADMM_multiple_primal} uses local residuals
$B_j v_j-s_j^{(n)}$ and intermediate multipliers
$\lambda^{(n+\frac{j}{J})}$, rather than the full constraint residual in each block subproblem.
The final step reconstructs the primal variables associated with $\lambda^{(n+1)}$ and does not affect the subsequent updates.

We next state the dualization relation.

\begin{theorem}
\label{Thm:ADMM_multiple}
In~\eqref{ADMM_multiple}, suppose that the functions $H_j$ defined in~\eqref{H_j} are convex.
Then the dualization-based ADMM~(\cref{Alg:ADMM_multiple_primal}) is a dualization of the Peaceman--Rachford splitting algorithm~\eqref{ADMM_multiple_PR}
for solving~\eqref{ADMM_multiple_dual_modified}.
More precisely, if
\begin{equation*}
\undertilde{s}^{(0)} = \undertilde{r}^{(0)}, \quad
\lambda^{(0)} = p^{(0)},
\end{equation*}
then we have
\begin{equation*}
\undertilde{s}^{(n)} = \undertilde{r}^{(n)}, \quad
\lambda^{(n)} = p^{(n)},
\quad
-B_j^*p^{(n)}\in\partial F_j(u_j^{(n)}),
\quad 1\leq j\leq J,\ n\geq 1.
\end{equation*}
\end{theorem}

\begin{proof}
Assume that $\undertilde{s}^{(n)}=\undertilde r^{(n)}$ and
$\lambda^{(n)}=p^{(n)}$.
We prove by induction over $j$ that
\[
\lambda^{(n+\frac{j}{J})}=p^{(n+\frac{j}{J})},
\quad 0\leq j\leq J.
\]
The case $j=0$ is immediate. Suppose the identity holds for $j-1$.
Using~\eqref{H_j}, the $j$th subproblem in~\eqref{ADMM_multiple_PR} is, up to constants independent of $q$,
\[
\min_{q\in W^*}
\left\{
F_j^*(-B_j^*q)
+
\frac{1}{2\gamma_j\beta}\|q\|^2
+
\left(
q,
s_j^{(n)}
-
\frac{2}{\beta}\lambda^{(n+\frac{j-1}{J})}
\right)
\right\}.
\]
Its Fenchel--Rockafellar dual is equivalent to the $\hat u_j^{(n+1)}$-subproblem in~\cref{Alg:ADMM_multiple_primal}, and the first primal--dual relation gives
\[
B_j\hat u_j^{(n+1)}
=
s_j^{(n)}
-
\frac{2}{\beta}p^{(n+\frac{j-1}{J})}
+
\frac{1}{\gamma_j\beta}p^{(n+\frac{j}{J})}.
\]
Substituting this identity into the update of
$\lambda^{(n+\frac{j}{J})}$ gives
\[
\lambda^{(n+\frac{j}{J})}=p^{(n+\frac{j}{J})}.
\]
The update of $s_j^{(n+1)}$ then gives
\[
s_j^{(n+1)}
=
s_j^{(n)}
+
\frac{2}{\beta}
(
p^{(n+\frac{j}{J})}
-
p^{(n+\frac{j-1}{J})}
)
=
r_j^{(n+1)}.
\]
Thus $\undertilde{s}^{(n+1)}=\undertilde r^{(n+1)}$ and
$\lambda^{(n+1)}=p^{(n+1)}$.
The final step of \cref{Alg:ADMM_multiple_primal} therefore gives
$-B_j^*p^{(n+1)}\in\partial F_j(u_j^{(n+1)})$ for $1\leq j\leq J$.
The result follows by induction.
\end{proof}

As a straightforward corollary of~\cref{Thm:multiple_convex_PR,Thm:ADMM_multiple}, we obtain the following relation between SSC and dualization-based ADMM.

\begin{corollary}
\label{Cor:ADMM_multiple}
In~\eqref{ADMM_multiple}, suppose that the functions $H_j$ defined in~\eqref{H_j} are convex.
Let $\{ ( \undertilde{s}^{(n)}, \lambda^{(n)} ) \}$ and $\{ \undertilde{w}^{(n)} \}$ be the sequences generated by the dualization-based ADMM~(\cref{Alg:ADMM_multiple_primal}) and SSC with exact local problems applied to the setting
\begin{equation*}
V \leftarrow W^J, \quad
E(\undertilde{w}) \leftarrow
\frac{\beta}{4}
\left\|
    \sum_{j=1}^J w_j - g
\right\|^2
+
\sum_{j=1}^J H_j^*(-w_j),
\end{equation*}
respectively, where $\undertilde{w}=(w_1,\dots,w_J)\in W^J$.
If
\begin{equation*}
    \lambda^{(0)}
    =
    \frac{\beta}{2}
    \left(
        \sum_{j=1}^J w_j^{(0)} - g
    \right),
    \quad
    s_j^{(0)} = w_j^{(0)},
    \quad 1\leq j\leq J,
\end{equation*}
then we have
\begin{equation*}
    \lambda^{(n)}
    =
    \frac{\beta}{2}
    \left(
        \sum_{j=1}^J w_j^{(n)} - g
    \right),
    \quad
    s_j^{(n)} = w_j^{(n)},
    \quad 1\leq j\leq J,\quad n\geq 1.
\end{equation*}
\end{corollary}

\begin{remark}
\label{Rem:ADMM_one_strongly_convex_block}
The requirement that both
\(F_{J-1}^*\circ(-B_{J-1}^*)\) and
\(F_J^*\circ(-B_J^*)\) be strongly convex is only a convenient sufficient condition.
It can be replaced by strong convexity of a single such term: one may split that term into two parts and apply the same construction to the resulting \(J+1\)-term dual problem, with \(\beta>0\) chosen sufficiently large.
\end{remark}

\section{Concluding remarks}
\label{Sec:Conclusion}
We have established connections among three broad classes of convex optimization algorithms based on the principle of \emph{divide, conquer, and combine}: subspace correction methods, operator splitting methods~(including alternating projection methods), and multiplier methods.
Despite their apparent differences, the operator splitting, alternating projection, and multiplier algorithms examined in the preceding sections can all be traced to subspace correction through convex duality and equivalent block formulations.

The relationships among the principal algorithms are summarized in \cref{Fig:algorithm_relationships}.
Starting from PSC and SSC, the expanded block methods are obtained through equivalent block formulations, while the splitting methods are related through extension and dualization.
Further dualization leads to the multiplier methods.
Thus, every algorithm in the figure is connected to a subspace correction method through extension, dualization, or equivalence.

\begin{figure}[t]
\centering
\resizebox{\textwidth}{!}{%
\begin{tikzpicture}[
algorithm node/.style={
    draw,
    rounded corners=2pt,
    align=center,
    minimum width=3.15cm,
    minimum height=0.85cm,
    inner sep=3pt,
    font=\small
},
region/.style={
    draw=black!60,
    rounded corners=5pt,
    line width=0.7pt
},
subspace region/.style={region,fill=blue!7},
splitting region/.style={region,fill=green!7},
admm region/.style={region,fill=orange!10},
extension/.style={-{Stealth[length=2.2mm]},draw=blue!70!black,line width=0.8pt},
dualization/.style={-{Stealth[length=2.2mm]},draw=red!75!black,line width=0.8pt,dashed},
equivalence/.style={draw=green!45!black,double=white,double distance=1.1pt,line width=0.5pt}
]

\begin{scope}[on background layer]
    \draw[subspace region] (-1.9,6.9) rectangle (14.2,10.6);
    \draw[splitting region] (-1.9,2.1) rectangle (14.2,6.55);
    \draw[admm region] (-1.9,-0.9) rectangle (14.2,1.8);
\end{scope}

\node[font=\bfseries] at (6.15,10.25) {Subspace correction methods};
\node[font=\bfseries] at (6.15,6.2) {Operator splitting methods};
\node[font=\bfseries] at (6.15,1.45) {ADMMs};

\node[algorithm node] (psc) at (0,8.55)
    {Parallel\\subspace correction\\\cref{Alg:MSC_linear_PSC,Alg:PSC}};
\node[algorithm node] (ebj) at (4.1,8.55)
    {Expanded\\Jacobi method\\\cref{Alg:MSC_linear_block_Jacobi,Alg:MSC_block_Jacobi}};
\node[algorithm node] (ssc) at (8.2,8.55)
    {Successive\\subspace correction\\\cref{Alg:MSC_linear_SSC,Alg:SSC}};
\node[algorithm node] (ebgs) at (12.3,8.55)
    {Expanded\\Gauss--Seidel method\\\cref{Alg:MSC_linear_block_GS,Alg:MSC_block_GS}};

\node[algorithm node] (pdr) at (0,5.1)
    {Parallel\\Douglas--Rachford\\splitting algorithm\\\cref{Alg:multiple_convex_parallel}};
\node[algorithm node] (pr) at (4.1,5.1)
    {Peaceman--Rachford\\splitting algorithm\\\cref{Alg:multiple_convex_PR}};
\node[algorithm node] (dr) at (8.2,5.1)
    {Douglas--Rachford\\splitting algorithm\\\cref{Alg:multiple_convex_DR}};
\node[algorithm node] (dykstra) at (4.1,3.25)
    {Dykstra algorithm\\\cref{Alg:Dykstra}};
\node[algorithm node] (vonneumann) at (8.2,3.25)
    {von Neumann algorithm\\\cref{Alg:MSC_von_Neumann}};
\node[algorithm node] (kaczmarz) at (12.3,3.25)
    {Kaczmarz algorithm\\\cref{Alg:Kaczmarz}};

\node[algorithm node] (pradmm) at (4.1,0.2)
    {Dualization-based\\ADMM\\\cref{Alg:ADMM_multiple_primal}};
\node[algorithm node] (admm) at (8.2,0.2)
    {Plain ADMM\\\cref{Alg:ADMM_multiple}};
\node[algorithm node] (sadmm) at (12.3,0.2)
    {Symmetrized ADMM\\\cref{Alg:ADMM_multiple_symmetrized}};

\draw[equivalence] (psc) -- (ebj);
\draw[equivalence] (ssc) -- (ebgs);

\draw[dualization] (psc) -- (pdr);
\draw[dualization,rounded corners=4pt]
    (ssc.south west) -- ++(0,-0.45) -- (2.8,7.25)
    -- (2.8,6.1) -- (pr.north west);
\draw[dualization] (ssc.south east) to[out=-90,in=90] (dr.north east);
\draw[extension] (dykstra) -- (pr);
\draw[extension] (vonneumann) -- (dykstra);
\draw[equivalence] (kaczmarz) -- (vonneumann);

\draw[dualization] (pr.west) -- ++(-0.45,0) |- (pradmm.west);
\draw[dualization] (dr.east) -- ++(0.45,0)
    |- ([yshift=-0.35cm]admm.south) -- (admm.south);
\node[font=\scriptsize,text=red!75!black,fill=green!7,inner sep=1pt]
    at (10.225,4.15) {$J=2$};
\draw[extension] (admm) -- (sadmm);

\end{tikzpicture}%
}
\caption{Relationships among the principal algorithms considered in the paper.
A directed edge points from a source algorithm to an algorithm obtained from it by extension or dualization, whereas an undirected edge denotes equivalence.
Solid blue, dashed red, and double green edges denote extension, dualization, and equivalence, respectively.
}
\label{Fig:algorithm_relationships}
\end{figure}

This paper suggests several interesting directions for future research.
First, the subspaces used in subspace correction methods may overlap, unlike the disjoint variable blocks typically used in conventional block methods.
Through the framework developed in this paper, this feature may be transferred to convex optimization to design efficient new algorithms.
Indeed, overlapping subspaces play an important role in scientific computing, particularly in domain decomposition and multigrid methods~\cite{Bramble:1993,TW:2005,XZ:1998}.

Second, the expanded system formulation allows subspace correction methods to be interpreted as block methods.
This viewpoint makes it possible to transfer tools from block coordinate descent to scientific computing, including randomized and importance-sampling strategies~\cite{FR:2015,FR:2016,RT:2016}, acceleration~\cite{CP:2015,LLX:2015,LX:2015}, and upper-bound minimization~\cite{HWRL:2017,RHL:2013}.
These tools may provide useful ingredients for randomized, accelerated, and nonlinear subspace correction methods.

Finally, refined and more general convergence theories are needed for subspace correction methods for convex optimization.
As discussed in \cref{Sec:MSC}, a sharp convergence theory for SSC that directly extends the linear theory is not yet available; see, e.g.,~\cite{BK:2012,CHW:2020,TE:1998} for existing convergence analyses.
Combined with the algorithmic connections presented in this paper, such results would support a common analysis of a broad range of convex optimization algorithms.

\bibliographystyle{siamplain}
\bibliography{refs_duality}

\end{document}